\documentclass[reqno]{amsart}
\usepackage{amscd}
\usepackage{amssymb}
\usepackage[dvips]{graphics}
\usepackage[margin=1in,includeheadfoot]{geometry}
\usepackage{url}
\usepackage{hyperref}
\usepackage{cite}
\hypersetup{hidelinks}

\def\beq{\begin{equation}}
\def\eeq{\end{equation}}
\def\ba{\begin{array}}
\def\ea{\end{array}}

\newtheorem{thm}{Theorem}[section]
\newtheorem{lm}[thm]{Lemma}
\newtheorem{prop}[thm]{Proposition}

\theoremstyle{definition}
\newtheorem{rem}[thm]{Remark}

\newtheorem{df}[thm]{Definition}
\newtheorem{ex}[thm]{Example}
\theoremstyle{remark}

\begin{document}
\pagestyle{plain}
\title{A new notion of Dirichlet boundary condition for elliptic equation of divergence form}

\author{Tianyu Guan\\
 \small{School of Mathematical Sciences, Shanghai Jiao Tong University}\\
 \small{Shanghai, China}\\\\
\small{ Lihe Wang}\\
 \small{Department of Mathematics, University of Iowa, Iowa City, IA, USA
 }\\\\
\small{Chunqin Zhou}\\
 \small{School of Mathematical Sciences, Shanghai Jiao Tong University}\\
 \small{Shanghai, China}}

\footnote{The third author was partially supported by NSFC Grant 12571223, 12031012 and STCSM Grant 24ZR1440700.}

\begin{abstract}
    In this paper, we consider Dirichlet boundary conditions for elliptic equations in divergence form on rough domains by introducing a truncation-based weak solution framework that circumvents global Sobolev trace obstacles. We establish well-posedness on bounded admissible domains for continuous boundary data and $L^2$
 sources, deriving boundary H\"older regularities that require only uniform ellipticity and an exterior measure-density condition. This approach eliminates the  uniform domain regularity and coefficient continuity imnposed in \cite{GWZ26}. Furthermore, on sufficiently flat Reifenberg domains, we obtain boundary $C^\alpha $ regularity for any $\alpha\in (0,1) $ using  natural energy estimates, thereby bypassing 
 $L^p$($p>\frac {n}{2}$) restrictions found in \cite{LZ26}. Crucially, we obtain global 
 regularity for H\"older boundary data ($\alpha>\frac 12$) a posteriori, without requiring small-BMO coefficient assumptions on the coefficients. 
\end{abstract}

\maketitle

\section{Introduction}

Let $\Omega\subset\mathbb R^n$ be a bounded domain, $n\geq 2$. We consider nonhomogeneous Dirichlet problem of the elliptic equation in divergence form 
\begin{equation}\label{1}
\left\{
\begin{array}{rcll}
-D_i\!\left(a_{ij}(x)D_j u\right)&=&f\qquad&\text{in }\Omega,\\
u&=&g&\text{on }\partial\Omega,
\end{array}
\right.
\end{equation}
where $A(x)=\{a_{ij}(x)\}_{i,j=1}^n$ is a real measurable matrix, and
repeated indices are summed. We assume throughout that, for some $\Lambda\geq1$,
\begin{equation}\label{ellipticity}
\Lambda^{-1}|\xi|^2\leq A(x)\xi\cdot\xi,
\qquad \|A\|_{L^\infty(\Omega)}\leq\Lambda
\end{equation}
for almost every $x\in \Omega $ and every $\xi\in\mathbb R^n$.
Here $|A|=(\sum_{i,j}|a_{ij}|^2)^{1/2}$ is the matrix norm used in
$\|A\|_{L^\infty}$ and in the small-oscillation condition below. 

For $g\in C(\partial \Omega)$, a classical problem concerns the solvability of the Dirichlet problem for the Laplace equation (i.e. f=0, A(x) is identity matrix) on arbitrary bounded domain $\Omega$.  In 1913, Lebesgue \cite{Leb13} constructed a bounded domain on which this is not solvable. In 1924, Wiener \cite{Wie24} proposed a necessary and sufficient condition for solvability, known as the Wiener criterion. This result implies that the domain must satisfy specific conditions to ensure the continuity of the solution up to the boundary. Subsequently, the Wiener criterion has been extended to linear and  quasilinear equations in divergence form;  see \cite{Kil94, Lit63, Lin85}  and the reference therein. 

On the other hand, on a Lipschitz domain, the boundary condition in \eqref{1} possesses the standard variational interpretation: one assumes $u\in H^1(\Omega)$ and prescribes the trace
of $u$ on $\partial\Omega$. However, on  rough domains, a canonical trace on the topological boundary may not exist. Consequently, it is unsatisfactory to impose  the global condition $u\in H^1(\Omega)$ in the
definition when global energy regularity is itself a conclusion to
be proved. Furthermore, for divergence-form equations, any replacement for the trace
condition must preserve the boundary energy inequalities upon which the
regularity theory relies.

Notivated by these consideration, in our former work \cite{GWZ26}, we introduced a definition of weak solutions for the Poisson problem on uniform domains. In this framework,  boundary data are encoded via the zero extensions of the upper and lower
boundary truncations, rather than through Sobolev traces. This definition requires only $H^1_{\mathrm{loc}}(\Omega)$ regularity in the interior while retaining the one-sided
boundary energy inequalities. However, the boundary
$C^\alpha$ theorem in \cite{GWZ26} assumes that the domain is both
$L$-uniform and $\theta$-admissible,  and it addresses ony the Laplace operator. Uniform
domains possess substantial interior geometry and, in particular, are Sobolev
extension domains; see \cite{Jon81}.

In this paper, we consider the more general equation \eqref{1},
employing the same solution definition as in \cite{GWZ26}. We first
address the question of well-posedness left open in \cite{GWZ26}.
For $g\in C(\partial\Omega)$ and $f\in L^2(\Omega)$, we establish
existence on every bounded domain and uniqueness on bounded
$\theta$-admissible domains.

We also consider  the boundary H\"older estimates. The primary advancement of the present paper is demonstrating that neither uniformity nor an extension-domain hypothesis is required for the basic boundary H\"older estimate. We retain only the exterior measure-density condition \eqref{admissible}. Moreover, the Laplacian in \cite{GWZ26} is replaced by a general divergence-form operator with coefficients that are measurable, bounded, and uniformly elliptic. No continuity, VMO, small-BMO, or pointwise closeness-to-a-constant assumptions on $A(x)$ are imposed in this part. Consequently, at the $C^\alpha$ level, the admissible-domain result presented here weakens both the geometric and coefficient hypotheses of our former work. The trade-off is that the available exponent is restricted to $0<\alpha<\beta(n,\Lambda,\theta)$.

Boundary H\"older continuity under an exterior cone or exterior
measure-density condition is classical; see, for example, \cite{GT01}.
Lemenant and Sire \cite{LS13} studied the Poisson equation on
Reifenberg-flat domains with zero boundary data, while Lian, Zhang, Li, and Hong \cite{LZLH20}
derived boundary H\"older estimates under a quantitative exterior geometric
condition. Additionally, Byun and Wang  \cite{BW04} established global $W^{1,p}$ estimates for small-BMO coefficients on
Reifenberg domains.

For sufficiently flat Reifenberg domains, Lian and Zhang \cite{LZ26}
recently obtained the boundary  H\"older estimate with arbitrary prescribed exponents $0<\alpha<1$ for several
classes of elliptic equations. Their theorem for divergence-form equations assumes
$u\in W^{1,2}(\Omega)\cap L^\infty(\Omega)$, interprets the boundary data in
the pointwise sense such that one can apply the maximum-principle, requires the coefficients to be
pointwise close to a constant matrix, and assumes
$f\in C^{-2,\alpha}_p$ for some $p>n/2$. In contrast, Our result on sufficiently flat Reifenberg domains provides an
$L^2$-Campanato estimate for weak solutions in the sense of Definition
\ref{uni}. It measures the right-hand side at the natural energy scale,
\[
 \|f\|_{L^2(\Omega_r(x_0))}
 \leq C r^{n/2-2+\alpha},
\]
and thereby elimination the need for an auxiliary exponent $p>n/2$. In dimensions $n\geq4$
this permits integrability at or below the threshold excluded in \cite{LZ26}; in all dimensions, it avoids imposing an  priori $L^\infty$ bound for the solution and operates directly with the energy norm. Furthermore, small mean oscillation replaces the requirement of pointwise closeness of $A$ to a fixed matrix.

%The restriction on the exponent in the first theorem is unavoidable at this
%level: arbitrarily high H\"older exponents cannot in general follow from
%measurable uniform ellipticity alone; see \cite{Mey63}.

Finally, on sufficiently flat Reifenberg domains, we obtain global
$H^1$ regularity for H\"older boundary data with exponent greater than
$1/2$. This lower threshold is sharp, as shown by the
classical Hadamard-type counterexample \cite{Had06,BDL26}:
Example \ref{ex:sharp-w1p}, with $p=2$, gives $C^{1/2}$ boundary
data on a smooth bounded domain in $\mathbb R^n$ whose harmonic extension has infinite Dirichlet
energy. The geometric boundary-layer estimate gives a continuous
$H^1$ lifting of the boundary data. We verify that the corresponding
variational solution satisfies our boundary truncation conditions and
then apply uniqueness. This argument requires only measurable uniformly
elliptic coefficients and $f\in L^2(\Omega)$; no small-BMO or boundary
Morrey condition is needed for the global $H^1$ regularity. Consequently  global $H^1$ regularity is derived a posteriori rather than assumed in the
definition.

The notation, geometric assumptions, and pointwise regularity spaces
are collected in Section \ref{sec:preliminaries}. For the coefficient
matrix, we use the following averages and small-oscillation condition.

For a ball $B_r(x)$, set
\[
(A)_{B_r(x)}:=\frac1{|B_r|}\int_{B_r(x)}A(y)\,dy.
\]
More generally, if $E\subset\mathbb R^n$ is measurable and $|E|>0$, set
\[
(A)_E:=\frac1{|E|}\int_EA(y)\,dy.
\]
The coefficient matrix $A$ is called $(\delta,R)$-vanishing if
\begin{equation}\label{small-bmo}
\sup_{0<r\leq R}\sup_{x\in\Omega}
\left(
\frac1{|B_r|}\int_{\Omega_r(x)}
\big|A(y)-(A)_{\Omega_r(x)}\big|^2\,dy
\right)^{1/2}
\leq\delta.
\end{equation}

We next state the definition of weak solutions used throughout the paper.

\begin{df}\label{uni}
Let $\Omega\subset\mathbb R^n$ be a bounded domain,  $A$ satisfies
\eqref{ellipticity}, $f\in L^2(\Omega)$, and $g\in L^\infty_{\mathrm{loc}}(\partial\Omega)$. A function $u$ is a weak
solution of \eqref{1} if the following conditions hold:
\begin{enumerate}
\renewcommand{\labelenumi}{(\roman{enumi})}
\item $u\in H^1_{\mathrm{loc}}(\Omega)$ and
\begin{equation}\label{eq1}
\int_\Omega a_{ij}(x)D_juD_i\varphi\,dx
=
\int_\Omega f\varphi\,dx
\qquad
\forall\,\varphi\in C_0^\infty(\Omega).
\end{equation}

\item For every $x_0\in\partial\Omega$ and
$0<r\leq\operatorname{diam}(\Omega)$, the zero extension in $\Omega^c$ of
\[
\left(u-\sup_{\partial_w\Omega_r(x_0)}g\right)^+
\]
belongs to $H^1_{\mathrm{loc}}(B_r(x_0))$ and satisfies
\begin{equation}\label{eq2}
\int_{B_r(x_0)}a_{ij}(x)D_j
\left[
\left(u-\sup_{\partial_w\Omega_r(x_0)}g\right)^+
\right]D_i\varphi\,dx
\leq
\int_{B_r(x_0)}f^+\varphi\,dx
\end{equation}
for every nonnegative $\varphi\in C_0^\infty(B_r(x_0))$.

\item For every $x_0\in\partial\Omega$ and
$0<r\leq\operatorname{diam}(\Omega)$, the zero extension in $\Omega^c$ of
\[
\left(u-\inf_{\partial_w\Omega_r(x_0)}g\right)^-
\]
belongs to $H^1_{\mathrm{loc}}(B_r(x_0))$ and satisfies
\begin{equation}\label{eq3}
\int_{B_r(x_0)}a_{ij}(x)D_j
\left[
\left(u-\inf_{\partial_w\Omega_r(x_0)}g\right)^-
\right]D_i\varphi\,dx
\leq
\int_{B_r(x_0)}f^-\varphi\,dx
\end{equation}
for every nonnegative $\varphi\in C_0^\infty(B_r(x_0))$.
\end{enumerate}
\end{df}
In the above definition, we extend $f(x)$ to zero in $\Omega ^{c}$.

We separate existence on arbitrary bounded domains from
uniqueness under the exterior measure-density condition.

\begin{thm}[Existence on bounded domains]\label{thm:existence}
Let $\Omega$ be a bounded domain, let $A$ satisfy
\eqref{ellipticity}, and let $f\in L^2(\Omega)$ and
$g\in C(\partial\Omega)$. Then problem \eqref{1} has a weak solution
$u\in L^2(\Omega)$ in the sense of Definition \ref{uni}, with
\begin{equation}\label{eq:existence-L2}
\|u\|_{L^2(\Omega)}
\leq C\|f\|_{L^2(\Omega)}
+|\Omega|^{1/2}\|g\|_{L^\infty(\partial\Omega)}.
\end{equation}
Here $C$ depends only on $n,\Lambda$, and
$\operatorname{diam}(\Omega)$. 
\end{thm}

\begin{thm}[Uniqueness on admissible domains]\label{thm:uniqueness}
Let $\Omega$ be a bounded $\theta$-admissible domain, let $A$
satisfy \eqref{ellipticity}, and let $f\in L^2(\Omega)$ and
$g\in C(\partial\Omega)$. Then the solution of \eqref{1} in the
sense of Definition \ref{uni} is unique.
\end{thm}

For existence, we approximate the continuous boundary data by smooth
functions and solve the variational problems on the fixed domain.
Comparison makes the approximating solutions Cauchy in $L^2(\Omega)$,
and the boundary energy estimate passes the truncation conditions to
the limit. 
For uniqueness, a weak solution is not assumed a priori to lie in
$H^1(\Omega)$. We construct upper and lower Sobolev approximations
near the boundary, prove a boundary-layer Poincar\'e inequality using
\eqref{admissible}, and apply a weak maximum principle.

Our first regularity theorem is the following.

\begin{thm}[Boundary pointwise $C^\alpha$ estimate on admissible domains]
\label{thm:boundary-holder-harmonic-replacement}
Let $\Omega$ be a bounded $\theta$-admissible domain and $0\in\partial\Omega$. Suppose $A$ satisfy \eqref{ellipticity}, and $u$ is a
weak solution of \eqref{1}. Then there exists
\[
\beta=\beta(n,\Lambda,\theta)\in(0,1)
\]
such that, for every $0<\alpha<\beta$, if
$f$ is $C^{-2,\alpha}$ at $0$ in the $L^2$ sense,
and $g$ is $C^\alpha$ at $0$, then $u$ is $C^\alpha$ at $0$ in the
$L^2$ sense, and
\[
[u]_{\mathcal L^{2,\alpha}(0)}
\leq C\left[
\left(
\frac1{|B_1|}\int_{\Omega_1}|u|^2\,dx
\right)^{1/2}
+\|g\|_{C^\alpha(0)}
+[f]_{\mathcal L^{2,-2+\alpha}(0)}
\right],
\]
where $C=C(n,\Lambda,\theta,\alpha)$.
\end{thm}

The proof combines the boundary energy estimate with the zero-extension
boundary De Giorgi--Moser decay for $A$-harmonic replacements. The one-sided
comparison with the original truncations is precisely that provided by Definition \ref{uni}.

The following theorem establishes a stronger, perturbative estimate on Reifenberg 
domains with small coefficient oscillation.

\begin{thm}[Boundary pointwise $C^\alpha$ estimate on  sufficiently flat Reifenberg domains]\label{thm1}
For any $0<\alpha<1$, there is a small $\delta^*=\delta^*(n,\alpha,\Lambda)>0$ so that for any $0<\delta<\delta^*$,  for all $A$ with A uniform elliptic(see \eqref{ellipticity}) and $(\delta,1)$-vanishing, for all $\Omega$ with
$(\delta,1)$-Reifenberg flat at $0\in\partial\Omega$, if $f$ is $C^{-2,\alpha}$ at $0$ in the $L^2$ sense,
and $g$ is $C^\alpha$ at $0$, then the
weak solution $u$ of \eqref{1} is $C^\alpha$ at $0$ in the
$L^2$ sense, with the estimate
\[
[u]_{\mathcal L^{2,\alpha}(0)}
\leq
C\left(
\frac1{|B_1|}\int_{\Omega_1}u^2\,dx
\right)^{1/2}
+
C\left(
[f]_{\mathcal L^{2,-2+\alpha}(0)}
+\|g\|_{C^\alpha(0)}
\right),
\]
where $C$ depends only on $n$, $\Lambda$, and $\alpha$.
\end{thm}

%The proof is an improvement-of-flatness iteration based on $L^2$
%approximation by a constant-coefficient solution in a half ball.

We now proceed to the recovery of  global $H^1(\Omega)$ regularity. Unlike the prescribed boundary H\"older exponent in Theorem \ref{thm1},
this conclusion does not require a small-BMO assumption on $A$.

\begin{thm}[Global $H^1(\Omega)$ regularity]\label{thm2}
Let $1/2<\alpha<1$ and $\Omega$ be a bounded
$(\delta,1)$-Reifenberg flat domain. Suppose that $A$ satisfies
\eqref{ellipticity}, $f\in L^2(\Omega)$, and
$g\in C^\alpha(\partial\Omega)$. Let $u$ be a weak solution of
\eqref{1} in the sense of Definition \ref{uni}. If
\[
0<\delta<\min\{\delta_0,\delta^{**}\},\qquad
\delta^{**}:=\sqrt{\frac{2\alpha-1}{C_1}},
\]
then $u\in H^1(\Omega)$ and
\[
\|u\|_{H^1(\Omega)}
\leq C\left(\|f\|_{L^2(\Omega)}
+\|g\|_{C^\alpha(\partial\Omega)}\right).
\]
Here $C_1=C_1(n)$ and $\delta_0=\delta_0(n)>0$ are the geometric
constants in Lemma \ref{rei}, with $\delta_0$ decreased if necessary
so that $\Omega$ satisfies the exterior measure-density condition
\eqref{admissible}. The constant $C$ depends only on $n$, $\Lambda$,
$\alpha$, $\delta$, and $\Omega$.
\end{thm}

The strict lower threshold $\alpha>1/2$ is sharp as mentioned above. The threshold $2\alpha>1+C_1\delta^{2}$ arises from the energy of a
Whitney extension of $g$. The gradient of this extension  is bounded by
$C\|g\|_{C^\alpha(\partial\Omega)}
\operatorname{dist}(x,\partial\Omega)^{\alpha-1}$, and the quantitative
size of boundary layers ensures that  this bound is square integrable. Theorem
\ref{thm:uniqueness} subsequently identifies the given weak solution with the
variational solution corresponding to this extension as boundary data. In
particular, this argument does not require  an a priori global $L^2$ bound on $u$.

Combining this construction with the global estimate of Byun and
Wang \cite[Theorem 1.5]{BW04} yields  global $W^{1,p}$ regularity.
Let $2\leq p<\infty$, $f\in L^p(\Omega)$, and
$g\in C^\alpha(\partial\Omega)$, where $0<\alpha<1$ and
$p(1-\alpha)<1-C_1\delta^{2}$. Under sufficiently small Reifenberg
flatness and the small-BMO assumption, depending on $p$, we obtain
\[
\|u\|_{W^{1,p}(\Omega)}
\leq C\left(\|f\|_{L^p(\Omega)}
+\|g\|_{C^\alpha(\partial\Omega)}\right).
\]
The boundary regularity of $g$ provides the necessary Whitney extension
in the proof. The precise statement and proof are detailed in Remark
\ref{rem:global-w1p}.
The lower threshold $\alpha>1-1/p$ is also sharp.

Our approach is significatly influenced by the energy method for equations with
rough coefficients on Reifenberg domains developed by Byun and Wang
\cite{BW04} as well as  by Caffarelli's compactness and perturbative methods
\cite{Caf88,Caf89} and their subsequent application to pointwise boundary regularity \cite{LiZ18,LZ20}. In essence, this  paper demonstrates that the truncation-based weak formulation from \cite{GWZ26} supports both a
nonperturbative boundary theory on admissible domains and a perturbative
boundary  theory on Reifenberg flat domains, while also establishing well-posedness and global $H^1$
regularity through a continuous energy lifting and uniqueness.

The paper is organized as follows. Section \ref{sec:preliminaries}
collects the domain geometry, boundary H\"older estimates, pointwise
norms, the Campanato embedding theorem, and the basic Whitney decomposition. Section \ref{sec:existence-uniqueness} first proves
Theorem \ref{thm:existence} on bounded domains, then proves Theorem
\ref{thm:uniqueness} on admissible domains. Section
\ref{sec:boundary-regularity} proves Theorem
\ref{thm:boundary-holder-harmonic-replacement}, using the boundary
energy estimate from Section \ref{sec:existence-uniqueness}. Section \ref{sec:reifenberg-regularity} addresses the
Reifenberg-domain setting: its first subsection proves Theorem \ref{thm1}
under the small-BMO assumption, and its second subsection proves
Theorem \ref{thm2} using only measurable uniformly elliptic coefficients.
The higher-integrability consequence under the small-BMO
assumption and the endpoint counterexample are recorded at the end of that subsection.

\section{Preliminaries}\label{sec:preliminaries}
In this section we collect the notation and standard facts used below.
We recall the geometric assumptions on the domain, the boundary H\"older
estimate, the pointwise norms and Campanato embedding theorem, and the
basic Whitney decomposition.

For $x_0\in\mathbb R^n$ and $r>0$, write
\[
 B_r(x_0)=\{x\in\mathbb R^n:|x-x_0|<r\},\qquad B_r=B_r(0),
\]
and
\[
 \Omega_r(x_0)=\Omega\cap B_r(x_0),\qquad \Omega_r=\Omega\cap B_r.
\]
For $x_0\in\partial\Omega$, the boundary portion is denoted by
\[
 \partial_w\Omega_r(x_0)=\partial\Omega\cap B_r(x_0).
\]
We omit the center when it is the origin, and write
$d(x,\partial\Omega)=\operatorname{dist}(x,\partial\Omega)$.
For $r>0$, the interior boundary layer is denoted by
\[
 \Omega^r:=\{x\in\Omega:d(x,\partial\Omega)<r\}.
\]
The right-hand side $f$ is extended by zero to $\Omega^c$ whenever it
occurs in an integral over a ball. Boundary truncations are also extended
by zero when considered on a ball; the extension is taken after truncation.
The letter $C$ denotes a positive
constant that may change from line to line; its dependence is specified
in the relevant statement or estimate.

\subsection{Geometric assumptions}
We first recall the exterior measure density condition.

\begin{df}[$\theta$-admissible domain]\label{df:admissible-domain}
Let $0<\theta<1$. A bounded domain $\Omega\subset\mathbb R^n$ is
called $\theta$-admissible if
\begin{equation}\label{admissible}
 |\Omega^c\cap B_r(x_0)|\geq\theta|B_r|
\end{equation}
for every $x_0\in\partial\Omega$ and every
$0<r<\operatorname{diam}(\Omega)$.
\end{df}
This is the geometric assumption used in the uniqueness
theorem and in the first boundary H\"older theorem.

For the stronger boundary estimate, we use Reifenberg flatness.

\begin{df}[Reifenberg flat domain]\label{df:reifenberg-domain}
Let $\delta,R>0$. A domain $\Omega\subset\mathbb R^n$ is
$(\delta,R)$-Reifenberg flat if, for every $x_0\in\partial\Omega$
and $0<r\leq R$, there is a coordinate system obtained by translating
$x_0$ to the origin and rotating the coordinates such that
\[
 B_r\cap\{x_n>\delta r\}
 \subset\Omega\cap B_r
 \subset B_r\cap\{x_n>-\delta r\}.
\]
\end{df}
As in \cite{BW04}, we also use the corresponding normalized local-coordinate
form
\[
B_r^+:=B_r\cap\{x_n>0\}
\subset\Omega\cap B_r
\subset B_r\cap\{x_n>-\delta r\}.
\]

\subsection{Boundary H\"older estimates}
The exterior density condition gives boundary H\"older decay for
solutions vanishing on the boundary in the Sobolev sense. We recall
the standard estimate; see \cite[Theorem 8.29]{GT01}.

\begin{lm}[Boundary H\"older estimate]\label{lm:preliminary-boundary-decay}
Let $\Omega$ be a $\theta$-admissible domain, let
$x_0\in\partial\Omega$, and let $0<R\leq\operatorname{diam}(\Omega)$.
Suppose that $A$ satisfies \eqref{ellipticity}, that
$h\in H^1(\Omega_R(x_0))$ satisfies
\[
 -D_i\!\left(a_{ij}(x)D_jh\right)=0
 \qquad\text{in }\Omega_R(x_0),
\]
and that its zero extension belongs to $H^1(B_R(x_0))$.
Then there exist $\beta=\beta(n,\Lambda,\theta)\in(0,1)$ and
$C=C(n,\Lambda,\theta)$ such that
\[
 |h(x)|\leq C\left(\frac{|x-x_0|}{R}\right)^\beta
 \left(\frac1{|B_R|}\int_{\Omega_R(x_0)}|h|^2\,dx\right)^{1/2}
 \qquad\bigl(x\in\Omega_{R/4}(x_0)\bigr),
\]
where $h$ is taken to be its continuous representative in the domain.
\end{lm}
In particular, for $0<r\leq R/4$,
\[
 \left(\frac1{|B_r|}\int_{\Omega_r(x_0)}|h|^2\,dx\right)^{1/2}
 \leq C\left(\frac rR\right)^\beta
 \left(\frac1{|B_R|}\int_{\Omega_R(x_0)}|h|^2\,dx\right)^{1/2}.
\]
We use this integral form in the boundary regularity proof.

\subsection{Pointwise norms and Campanato embedding}
We next recall the pointwise norms appearing in the main results.
Let $x_0\in\partial\Omega$, $r_0>0$, and $1\leq q<\infty$. For $\gamma<0$,
define
\[
 [f]_{\mathcal L^{q,\gamma}(x_0)}
 :=\sup_{0<r<r_0}r^{-\gamma}
 \left(\frac1{|B_r|}\int_{\Omega_r(x_0)}|f|^q\,dx\right)^{1/q}.
\]
In particular, for $0<\alpha<1$, we say that $f$ is $C^{-2,\alpha}$
at $x_0$ in the $L^2$ sense if
$[f]_{\mathcal L^{2,-2+\alpha}(x_0)}<\infty$.

For $0<\alpha<1$, we say that $u$ is $C^\alpha$ at $x_0$ in the
$L^q$ sense if there exists a constant $c\in\mathbb R$ and $r_0>0$ such that
\[
 \sup_{0<r<r_0}r^{-\alpha}
 \left(\frac1{|B_r|}\int_{\Omega_r(x_0)}|u-c|^q\,dx\right)^{1/q}
 <\infty.
\]
The corresponding pointwise seminorm is
\[
 [u]_{\mathcal L^{q,\alpha}(x_0)}
 :=\inf_{c\in\mathbb R}\sup_{0<r<r_0}r^{-\alpha}
 \left(\frac1{|B_r|}\int_{\Omega_r(x_0)}|u-c|^q\,dx\right)^{1/q}.
\]
In the boundary estimates below the approximating constant is $g(x_0)$.
For boundary data, we use
\[
 [g]_{C^\alpha(x_0)}
 :=\sup_{y\in\partial_w\Omega_1(x_0)\setminus\{x_0\}}
 \frac{|g(y)-g(x_0)|}{|y-x_0|^\alpha},
\]
and
\[
 \|g\|_{C^\alpha(x_0)}:=|g(x_0)|+[g]_{C^\alpha(x_0)}.
\]
For the global boundary norm, we use
\[
\|g\|_{C^\alpha(\partial\Omega)}
:=\|g\|_{L^\infty(\partial\Omega)}
+\sup_{\substack{y,z\in\partial\Omega\\y\ne z}}
\frac{|g(y)-g(z)|}{|y-z|^\alpha}.
\]

We also recall the degree-zero case of the Campanato embedding theorem;
see \cite[pp.~182--184]{Camp63}.

\begin{thm}[Campanato embedding]\label{thm:campanato-embedding}
Let $E\subset\mathbb R^n$ be a bounded open set satisfying
\[
 |E\cap B_r(x)|\geq c_0r^n
 \qquad\bigl(x\in\overline E,\ 0<r<r_0\bigr)
\]
for some $c_0,r_0>0$. Let $0<\alpha<1$ and $u\in L^2(E)$.
Suppose that, for every $x\in\overline E$ and $0<r<r_0$, there
exists a constant $c_{x,r}\in\mathbb R$ such that
\[
 \left(\frac1{|E\cap B_r(x)|}
 \int_{E\cap B_r(x)}|u-c_{x,r}|^2\,dy\right)^{1/2}
 \leq Mr^\alpha,
\]
where $M$ is independent of $x$ and $r$. Then $u$ has a representative
in $C^\alpha(\overline E)$, and
\[
 \|u\|_{C^\alpha(\overline E)}
 \leq C\bigl(\|u\|_{L^2(E)}+M\bigr),
\]
where $C$ depends only on $n,\alpha,c_0$, and $r_0$.
\end{thm}

\subsection{Whitney decomposition}
We recall the Whitney decomposition and the covering properties used in
the global energy argument; see \cite[Chapter VI, Section 1]{Stein70} for
the decomposition and the bounded overlap of the dilated cubes.

\begin{lm}[Whitney decomposition]\label{wd}
Let $F$ be a nonempty closed set in $\mathbb{R}^n$ and $\Omega$ be its complement. Then there exists a sequence of cubes $Q_k$ such that
\begin{enumerate}
    \item $\Omega=\bigcup_{k=1}^{\infty}Q_k$
    \item  The $Q_k$ are mutually disjoint;
    \item $d_k\leq \operatorname{dist}(Q_k,F)\leq 4d_k$, where $d_k=\operatorname{diam} Q_k$
\end{enumerate}
\end{lm}
\begin{lm}\label{wd2}
    Let $Q_k$ be as in Lemma $\ref{wd}$ and $\widetilde Q_k=\frac{6}{5}Q_k$ be $\frac{6}{5}$-dilation of $Q_k$ with respect to its center. Then
    \begin{enumerate}
        \item $\Omega=\bigcup_{k=1}^{\infty}\widetilde Q_k$;
        \item Each point of $\Omega$ is contained in at most $12^n$ of the cubes $\widetilde Q_k$.
    \end{enumerate}
    \end{lm}
\begin{lm}\label{wd3}
\begin{align}
    \Omega_{\frac{r}{3}}\subset \bigcup_{\widetilde Q_k\subset\Omega_r}Q_k
    \quad\text{for }0<r\leq1.
\end{align}
\end{lm}

\section{Existence and uniqueness}\label{sec:existence-uniqueness}
We first record two boundary estimates that require no geometric
assumption. We then construct solutions on arbitrary bounded domains.
The uniqueness argument follows and uses the exterior measure-density
condition.

\begin{lm}[Boundary energy estimate]\label{lm1}
If $u$ is a weak solution of \eqref{1} in $\Omega$, then for every $\eta\in C_0^\infty(B_1)$,
\begin{align}
\int_{\Omega_1}\eta^2|\nabla (u-\sup_{\partial_w\Omega_1}g)^+|^2\,dx
&\leq \int_{\Omega_1}\eta^2f^2\,dx
+C\int_{\Omega_1}(|\nabla\eta|^2+\eta^2)|u|^2\,dx \notag\\
&\quad
+C\int_{\Omega_1}(|\nabla\eta|^2+\eta^2)\,dx\,
\|g\|_{L^\infty(\partial_w\Omega_1)}^2, \label{v1}\\
\int_{\Omega_1}\eta^2|\nabla (u-\inf_{\partial_w\Omega_1}g)^-|^2\,dx
&\leq \int_{\Omega_1}\eta^2f^2\,dx
+C\int_{\Omega_1}(|\nabla\eta|^2+\eta^2)|u|^2\,dx \notag\\
&\quad
+C\int_{\Omega_1}(|\nabla\eta|^2+\eta^2)\,dx\,
\|g\|_{L^\infty(\partial_w\Omega_1)}^2, \label{g1}
\end{align}
where $C$ depends only on $\Lambda$.
\end{lm}
\begin{proof}
Take $\eta^2(u-\sup_{\partial_w\Omega_1}g)^+$ as the test function. Then
\[
\begin{aligned}
&\int_{B_1}a_{ij}(x)D_j\bigl[(u-\sup_{\partial_w\Omega_1}g)^+\bigr]
D_i\bigl[\eta^2(u-\sup_{\partial_w\Omega_1}g)^+\bigr]\,dx\\
&\qquad\leq\int_{B_1}\eta^2(u-\sup_{\partial_w\Omega_1}g)^+f^+\,dx.
\end{aligned}
\]
By ellipticity and Cauchy's inequality,
\begin{align*}
&\frac{1}{\Lambda}\int_{B_1}\eta^2
 \bigl|\nabla(u-\sup_{\partial_w\Omega_1}g)^+\bigr|^2\,dx\\
&\quad\leq\int_{B_1}\eta^2(u-\sup_{\partial_w\Omega_1}g)^+f^+\,dx\\
&\qquad-2\int_{B_1}\eta(u-\sup_{\partial_w\Omega_1}g)^+
 a_{ij}(x)D_j\bigl[(u-\sup_{\partial_w\Omega_1}g)^+\bigr]D_i\eta\,dx\\
&\quad\leq\frac{1}{2\Lambda}\int_{B_1}\eta^2(f^+)^2\,dx
 +\frac{\Lambda}{2}\int_{B_1}\eta^2
 \bigl[(u-\sup_{\partial_w\Omega_1}g)^+\bigr]^2\,dx\\
&\qquad+\frac{1}{2\Lambda}\int_{B_1}\eta^2
 \bigl|\nabla(u-\sup_{\partial_w\Omega_1}g)^+\bigr|^2\,dx+C(\Lambda)\int_{B_1}|\nabla\eta|^2
 \bigl[(u-\sup_{\partial_w\Omega_1}g)^+\bigr]^2\,dx.
\end{align*}
Since $(u-\sup_{\partial_w\Omega_1}g)^+$ and $f^+$ vanish outside
$\Omega$ and
\[
\bigl[(u-\sup_{\partial_w\Omega_1}g)^+\bigr]^2
\leq C\bigl(|u|^2+\|g\|_{L^\infty(\partial_w\Omega_1)}^2\bigr)
\quad\text{in }\Omega_1,
\]
\eqref{v1} follows after absorption. 
The same argument applies to the second inequality.
\end{proof}

\begin{lm}[Further truncation]\label{lm:further-truncation}
Let $u$ be a weak solution in the sense of Definition \ref{uni}.
Fix $x_0\in\partial\Omega$ and
$0<r\leq\operatorname{diam}(\Omega)$. If
\[
 c>\sup_{\partial_w\Omega_r(x_0)}g,
\]
the zero extension in $\Omega^c$ of $(u-c)^+$ belongs to
$H^1_{\mathrm{loc}}(B_r(x_0))$ and satisfies \eqref{eq2} with
$\sup_{\partial_w\Omega_r(x_0)}g$ replaced by $c$.
If
\[
 c<\inf_{\partial_w\Omega_r(x_0)}g,
\]
the zero extension in $\Omega^c$ of $(u-c)^-$ belongs to
$H^1_{\mathrm{loc}}(B_r(x_0))$ and satisfies \eqref{eq3} with
$\inf_{\partial_w\Omega_r(x_0)}g$ replaced by $c$.
\end{lm}
\begin{proof}
We first consider $c>\sup_{\partial_w\Omega_r(x_0)}g$. The identity
\[
(u-c)^+=\left((u-\sup_{\partial_w\Omega_r(x_0)}g)^+
-\left(c-\sup_{\partial_w\Omega_r(x_0)}g\right)\right)^+
\]
holds in $\Omega_r(x_0)$ and for the zero extensions to $B_r(x_0)$. Definition
\ref{uni} therefore gives
\[
(u-c)^+\in H^1_{\mathrm{loc}}(B_r(x_0)),\qquad
\nabla(u-c)^+=\mathbf{1}_{\{(u-c)^+>0\}}
\nabla(u-\sup_{\partial_w\Omega_r(x_0)}g)^+.
\]

Fix $0\leq\varphi\in C_c^\infty(B_r(x_0))$ and $\varepsilon>0$.
Take $\varphi\min\{(u-c)^+/\varepsilon,1\}$ as the test function in \ref{eq2}, then
\[
\begin{aligned}
&\int_{B_r(x_0)}\min\{(u-c)^+/\varepsilon,1\}
  A\nabla(u-c)^+\cdot\nabla\varphi\,dx\\
&\quad+\frac1\varepsilon
  \int_{\Omega_r(x_0)\cap\{c<u<c+\varepsilon\}}
  \varphi A\nabla(u-c)^+\cdot\nabla(u-c)^+\,dx\\
&\qquad\leq\int_{B_r(x_0)}f^+\varphi
  \min\{(u-c)^+/\varepsilon,1\}\,dx.
\end{aligned}
\]
The second integral on the left is nonnegative by ellipticity.
Dropping it and letting $\varepsilon\downarrow0$, dominated convergence
and the fact that $\nabla(u-c)^+=0$ almost everywhere on
$\{(u-c)^+=0\}$ yield
\[
\int_{B_r(x_0)}A\nabla(u-c)^+\cdot\nabla\varphi\,dx
\leq\int_{\Omega_r(x_0)\cap\{u>c\}}f^+\varphi\,dx
\leq\int_{B_r(x_0)}f^+\varphi\,dx.
\]
This is the required inequality for the zero extension of $(u-c)^+$.

For $c<\inf_{\partial_w\Omega_r(x_0)}g$, use the identity
\[
(u-c)^-=\left((u-\inf_{\partial_w\Omega_r(x_0)}g)^-
-\left(\inf_{\partial_w\Omega_r(x_0)}g-c\right)\right)^+.
\]
The same calculation applied to \eqref{eq3} gives the result with
right-hand side $f^-$.
\end{proof}

\subsection{Existence on bounded domains}

\begin{proof}[Proof of Theorem \ref{thm:existence}]
By the Tietze extension theorem, extend $g$ continuously to $\mathbb R^n$ and multiply the extension by a smooth compactly supported cutoff equal to one near $\overline{\Omega}$ to obtain $G\in C_c(\mathbb R^n)$. Then choose
$G_j\in C_c^\infty(\mathbb R^n)$ converging uniformly to $G$.
Lax--Milgram and Poincar\'e's inequality give $u_j\in H^1(\Omega)$
with $u_j-G_j\in H_0^1(\Omega)$ satisfying
\begin{equation}\label{eq:linear-approximation-equation}
\int_\Omega A\nabla u_j\cdot\nabla\psi\,dx
=\int_\Omega f\psi\,dx
\qquad\forall\psi\in H_0^1(\Omega).
\end{equation}
We verify the boundary conditions in Definition \ref{uni} for $u_j$.
Fix $x_0\in\partial\Omega$, $0<r\leq\operatorname{diam}(\Omega)$, and
$\eta\in C_0^\infty(B_r(x_0))$. Choose $w_m\in C_0^\infty(\Omega)$
converging to $u_j-G_j$ in $H^1(\Omega)$. For $\varepsilon>0$,
\[
\eta\left(G_j+w_m-\sup_{\partial_w\Omega_r(x_0)}G_j-\varepsilon\right)^+
\]
has compact support in $\Omega$, since $w_m$ vanishes near the
boundary and continuity of $G_j$ makes the positive part vanish near
$\partial\Omega\cap\operatorname{supp}\eta$. Letting $m\to\infty$ and
then $\varepsilon\downarrow0$, and applying the same argument to
$-u_j,-G_j$, gives
\begin{equation}\label{eq:variational-local-truncation}
\eta\left(u_j-\sup_{\partial_w\Omega_r(x_0)}G_j\right)^+
\in H_0^1(\Omega),\qquad
\eta\left(u_j-\inf_{\partial_w\Omega_r(x_0)}G_j\right)^-
\in H_0^1(\Omega).
\end{equation}
Here we used continuity of truncation in $H^1(\Omega)$ and
closedness of $H_0^1(\Omega)$. Taking $\eta=1$ on smaller balls
shows that both truncations extend by zero to $H^1_{\mathrm{loc}}(B_r(x_0))$.
For $0\leq\varphi\in C_0^\infty(B_r(x_0))$, choose $\eta\geq0$
equal to one near $\operatorname{supp}\varphi$ in
\eqref{eq:variational-local-truncation}. Truncation then gives the
admissible test function
\[
\varphi\min\!\left\{1,\frac{\bigl(u_j-\sup_{\partial_w\Omega_r(x_0)}G_j\bigr)^+}{t}\right\},
\qquad t>0.
\]
By \eqref{eq:linear-approximation-equation},
\[
\begin{aligned}
&\int_{\Omega_r(x_0)}
\min\!\left\{1,\frac{\bigl(u_j-\sup_{\partial_w\Omega_r(x_0)}G_j\bigr)^+}{t}\right\}
A\nabla u_j\cdot\nabla\varphi\,dx\\
&\quad+
\underbrace{\frac1t
\int_{\Omega_r(x_0)\cap\{0<u_j-\sup_{\partial_w\Omega_r(x_0)}G_j<t\}}
\varphi A\nabla u_j\cdot\nabla u_j\,dx}_{\geq0}\\
&=\int_{\Omega_r(x_0)}f\varphi
\min\!\left\{1,\frac{\bigl(u_j-\sup_{\partial_w\Omega_r(x_0)}G_j\bigr)^+}{t}\right\}\,dx
\leq\int_{B_r(x_0)}f^+\varphi\,dx.
\end{aligned}
\]
Hence
\[
\begin{aligned}
&\int_{B_r(x_0)}
A\nabla\left(u_j-\sup_{\partial_w\Omega_r(x_0)}G_j\right)^+
\cdot\nabla\varphi\,dx\\
&=\lim_{t\downarrow0}\int_{\Omega_r(x_0)}
\min\!\left\{1,\frac{\bigl(u_j-\sup_{\partial_w\Omega_r(x_0)}G_j\bigr)^+}{t}\right\}
A\nabla u_j\cdot\nabla\varphi\,dx\\
&\leq\int_{B_r(x_0)}f^+\varphi\,dx.
\end{aligned}
\]
This proves \eqref{eq2} for $u_j$.
Applying this argument to $-u_j$ gives \eqref{eq3}, so $u_j$ satisfies
Definition \ref{uni} with boundary datum $G_j|_{\partial\Omega}$.

Fix $j,\ell$. Since $(u_j-G_j)-(u_\ell-G_\ell)\in H_0^1(\Omega)$,
choose $\psi_m\in C_0^\infty(\Omega)$ converging to this difference
in $H^1(\Omega)$. For each $\varepsilon>0$, the function
\[
\left(\psi_m+G_j-G_\ell-\sup_{\partial\Omega}|G_j-G_\ell|-\varepsilon\right)^+
\]
has compact support in $\Omega$. Indeed, $\psi_m$ vanishes near
$\partial\Omega$, and continuity of $G_j-G_\ell$ gives
$G_j-G_\ell<\sup_{\partial\Omega}|G_j-G_\ell|+\varepsilon$ in a
neighborhood of the compact boundary. This positive part thus belongs
to $H_0^1(\Omega)$ by mollification. As $m\to\infty$, it converges in
$H^1(\Omega)$ to
$(u_j-u_\ell-\sup_{\partial\Omega}|G_j-G_\ell|-\varepsilon)^+$
by continuity of the positive-part map. Letting $\varepsilon\downarrow0$
and using the closedness of $H_0^1(\Omega)$ in $H^1(\Omega)$ gives
\[
\left(u_j-u_\ell-\sup_{\partial\Omega}|G_j-G_\ell|\right)^+
\in H_0^1(\Omega).
\]
By \eqref{eq:linear-approximation-equation} and Poincar\'e's inequality, we obtain
\[
\begin{aligned}
&\left\|\left(u_j-u_\ell-\sup_{\partial\Omega}|G_j-G_\ell|\right)^+\right\|_{L^2(\Omega)}^2\\
&\quad\leq C\int_\Omega
\left|\nabla\left(u_j-u_\ell-\sup_{\partial\Omega}|G_j-G_\ell|\right)^+\right|^2\,dx\\
&\quad\leq C\int_\Omega A\nabla(u_j-u_\ell)\cdot
\nabla\left(u_j-u_\ell-\sup_{\partial\Omega}|G_j-G_\ell|\right)^+\,dx=0.
\end{aligned}
\]
Interchanging $j$ and $\ell$ gives

\begin{equation}\label{eq:linear-approximation-stability}
\|u_j-u_\ell\|_{L^\infty(\Omega)}
\leq\sup_{\partial\Omega}|G_j-G_\ell|.
\end{equation}
Let $u_0\in H_0^1(\Omega)$ be the variational solution with
right-hand side $f$ and zero boundary data. The same comparison gives
\[
|u_j-u_0|\leq\sup_{\partial\Omega}|G_j|
\quad\text{a.e. in }\Omega.
\]
Testing with $u_0$ in $u_{0}$'s equation and applying Poincar\'e's inequality give
\[
\begin{aligned}
\|u_j\|_{L^2(\Omega)}
&\leq C\|f\|_{L^2(\Omega)}
+|\Omega|^{1/2}\sup_{\partial\Omega}|G_j|,\\
\|u_j-u_\ell\|_{L^2(\Omega)}
&\leq|\Omega|^{1/2}\sup_{\partial\Omega}|G_j-G_\ell|
\longrightarrow0\qquad(j,\ell\to\infty).
\end{aligned}
\]
Hence $u_j\to u$ in $L^2(\Omega)$, and \eqref{eq:existence-L2} follows.

For $\eta\in C_0^\infty(\Omega)$, testing the difference equation
with $\eta^2(u_j-u_\ell)$ gives
\[
\int_\Omega\eta^2|\nabla u_j-\nabla u_\ell|^2\,dx
\leq C\int_\Omega|\nabla\eta|^2|u_j-u_\ell|^2\,dx.
\]
Thus $u_j\to u$ strongly in $H^1_{\mathrm{loc}}(\Omega)$, and
passing to the limit proves the interior equation.

Fix $x_0\in\partial\Omega$ and $0<r\leq\operatorname{diam}(\Omega)$.
Uniform convergence of $G_j$ to $g$ on $\partial\Omega$ gives
\[
\sup_{\partial_w\Omega_r(x_0)}G_j\longrightarrow\sup_{\partial_w\Omega_r(x_0)}g,
\qquad
\inf_{\partial_w\Omega_r(x_0)}G_j\longrightarrow\inf_{\partial_w\Omega_r(x_0)}g.
\]
Since $u_j\to u$ in $L^2(\Omega)$, the zero extensions of
$(u_j-\sup_{\partial_w\Omega_r(x_0)}G_j)^+$ converge strongly in
$L^2(B_r(x_0))$ to the zero extension of $(u-\sup_{\partial_w\Omega_r(x_0)}g)^+$.
The bounds on $\|u_j\|_{L^2(\Omega)}$ and $\sup_{\partial\Omega}|G_j|$
are uniform in $j$. Thus for every $0<s<r$, the boundary energy
estimate \eqref{v1}, after translation and rescaling, gives
\[
\sup_{j\geq1}\int_{\Omega_s(x_0)}
\left|\nabla\left(u_j-\sup_{\partial_w\Omega_r(x_0)}G_j\right)^+\right|^2\,dx
<\infty.
\]
Hence these zero extensions converge weakly in $H^1(B_s(x_0))$
to the same limit. For each $0\leq\varphi\in C_0^\infty(B_r(x_0))$,
choose $s<r$ with $\operatorname{supp}\varphi\subset B_s(x_0)$.
Since $A$ is bounded, passing to the limit in \eqref{eq2} for $u_j$ gives
\[
\int_{B_r(x_0)}A\nabla\left(u-\sup_{\partial_w\Omega_r(x_0)}g\right)^+
  \cdot\nabla\varphi\,dx\leq\int_{B_r(x_0)}f^+\varphi\,dx.
\]
Using \eqref{g1} and the convergence of the boundary infima in the
same argument shows that the zero extension of $(u-\inf_{\partial_w\Omega_r(x_0)}g)^-$
belongs to $H^1_{\mathrm{loc}}(B_r(x_0))$ and satisfies \eqref{eq3}.
Thus $u$
satisfies Definition \ref{uni}. 
\end{proof}

\subsection{Uniqueness on admissible domains}

We prove uniqueness using the boundary truncations in Definition
\ref{uni}. We first construct upper and lower $H^1$ approximations
near $\partial\Omega$.

\begin{prop}[Boundary approximation]\label{prop:boundary-approximation}
Let $g\in C(\partial\Omega)$ and let $u$ be a weak solution in the sense of
Definition \ref{uni}. For every $\varepsilon>0$, there exist
$u_\varepsilon^+,u_\varepsilon^-\in H^1(\Omega)$, an open neighborhood
$U_\varepsilon$ of $\partial\Omega$, functions
\[
\widetilde u_\varepsilon^+,\widetilde u_\varepsilon^-
\in H^1(U_\varepsilon),
\qquad
g_\varepsilon^+,g_\varepsilon^-\in C^\infty(U_\varepsilon),
\]
such that
\[
\widetilde u_\varepsilon^\pm=u_\varepsilon^\pm
\quad\text{a.e. in }U_\varepsilon\cap\Omega,
\qquad
\widetilde u_\varepsilon^\pm=g_\varepsilon^\pm
\quad\text{a.e. in }U_\varepsilon\cap\Omega^c,
\]
and
\[
u_\varepsilon^-\leq u\leq u_\varepsilon^+
\quad\text{a.e. in }U_\varepsilon\cap\Omega.
\]
Moreover,
\begin{equation}\label{eq:exterior-gap}
g_\varepsilon^+-g_\varepsilon^-=\varepsilon
\quad\text{in }U_\varepsilon,
\end{equation}
and
\begin{equation}\label{eq:boundary-approximation}
g-\varepsilon<g_\varepsilon^-<g<g_\varepsilon^+<g+\varepsilon
\quad\text{on }\partial\Omega.
\end{equation}
The neighborhood $U_\varepsilon$ and the functions
$g_\varepsilon^\pm$ depend only on $g$ and $\varepsilon$.
\end{prop}

\begin{proof}
Since $g$ is uniformly continuous on the compact set $\partial\Omega$, there
exist $x_1,\ldots,x_N\in\partial\Omega$ and
$0<2r_i<\operatorname{diam}(\Omega)$ such that
\[
\partial\Omega\subset\bigcup_{i=1}^N B_{r_i}(x_i)
\]
and
\begin{equation}\label{eq:g-oscillation}
|g(x)-g(x_i)|<\frac{\varepsilon}{4}
\quad\text{for every }x\in\partial_w\Omega_{2r_i}(x_i).
\end{equation}
Set
\[
c_i^+:=g(x_i)+\frac{\varepsilon}{2},
\qquad
c_i^-:=g(x_i)-\frac{\varepsilon}{2}.
\]
Then
\[
c_i^+>\sup_{\partial_w\Omega_{2r_i}(x_i)}g,
\qquad
c_i^-<\inf_{\partial_w\Omega_{2r_i}(x_i)}g.
\]

By Lemma \ref{lm:further-truncation}, the zero extensions of
$(u-c_i^+)^+$ and $(u-c_i^-)^-$ belong to
$H^1_{\mathrm{loc}}(B_{2r_i}(x_i))$. Their restrictions to the smaller
ball therefore satisfy
\[
(u-c_i^+)^+,(u-c_i^-)^-\in H^1(B_{r_i}(x_i)).
\]
In $B_{r_i}(x_i)$, set
\[
U_i^+:=c_i^++(u-c_i^+)^+,
\qquad
U_i^-:=c_i^--(u-c_i^-)^-.
\]
Then, a.e. in $B_{r_i}(x_i)\cap\Omega$,
\[
U_i^+=\max\{u,c_i^+\}\geq u,
\qquad
U_i^-=\min\{u,c_i^-\}\leq u,
\]
whereas, a.e. in $B_{r_i}(x_i)\cap\Omega^c$,
\[
U_i^+=c_i^+,\qquad U_i^-=c_i^-.
\]

Choose an open neighborhood $U$ of $\partial\Omega$ and a smooth partition
of unity $\{\eta_i\}_{i=1}^N$ in $U$, with $\operatorname{supp}\eta_i\Subset B_{r_i}(x_i)$, such that
\[
\sum_{i=1}^N\eta_i=1\quad\text{in }U.
\]
Define
\[
g_\varepsilon^+:=\sum_{i=1}^N\eta_i c_i^+,
\qquad
g_\varepsilon^-:=\sum_{i=1}^N\eta_i c_i^-,
\]
and
\[
\widetilde u_\varepsilon^+:=\sum_{i=1}^N\eta_iU_i^+,
\qquad
\widetilde u_\varepsilon^-:=\sum_{i=1}^N\eta_iU_i^-.
\]
Then $\widetilde u_\varepsilon^\pm\in H^1(U)$,
\[
\widetilde u_\varepsilon^-\leq u\leq
\widetilde u_\varepsilon^+\quad\text{a.e. in }U\cap\Omega,
\]
and
\[
\widetilde u_\varepsilon^\pm=g_\varepsilon^\pm
\quad\text{a.e. in }U\cap\Omega^c.
\]
Since $c_i^+-c_i^-=\varepsilon$, we have
\[
g_\varepsilon^+-g_\varepsilon^-=
\varepsilon\sum_{i=1}^N\eta_i=\varepsilon.
\]

If $x\in\partial\Omega$ and $\eta_i(x)\neq0$, then
\eqref{eq:g-oscillation} gives
\[
g(x)+\frac{\varepsilon}{4}<c_i^+
<g(x)+\frac{3\varepsilon}{4},
\]
and
\[
g(x)-\frac{3\varepsilon}{4}<c_i^-
<g(x)-\frac{\varepsilon}{4}.
\]
Taking convex combinations proves \eqref{eq:boundary-approximation}.

Finally, choose an open neighborhood $U_\varepsilon$ of $\partial\Omega$
with $\overline{U_\varepsilon}\subset U$ and
$\chi\in C_0^\infty(U)$ such that $\chi=1$ in $U_\varepsilon$. The
functions $\chi\widetilde u_\varepsilon^\pm$, extended by zero away from
$U$, belong to $H^1(\mathbb R^n)$. Their restrictions
$u_\varepsilon^\pm$ to $\Omega$ have the required properties.
The choices of $U_\varepsilon$ and $g_\varepsilon^\pm$ depend only on
$g$ and $\varepsilon$.
\end{proof}

\begin{lm}[Boundary-layer Poincar\'e inequality]\label{lm:boundary-poincare}
Let $\Omega$ be a bounded $\theta$-admissible domain. There exist $r_0>0$
and $C\geq1$, depending only on $n$, $\theta$, and $\Omega$, such that for
every $m\in H_0^1(\Omega)$ and every $0<r<r_0$,
\begin{equation}\label{eq:boundary-poincare}
\frac1{r^2}\int_{\Omega^r}|m|^2\,dx
\leq C\int_{\Omega^{2r}}|\nabla m|^2\,dx.
\end{equation}
Consequently,
\begin{equation}\label{eq:boundary-layer-vanishing}
\lim_{r\downarrow0}\frac1{r^2}\int_{\Omega^r}|m|^2\,dx=0.
\end{equation}
Moreover, if $z\in H^1(\mathbb R^n)$ and $z=0$ a.e. in $\Omega^c$, then
$z|_\Omega\in H_0^1(\Omega)$.
\end{lm}

\begin{proof}
Let $\widetilde m$ be the zero extension of $m$ to $\mathbb R^n$. Fix $r>0$
sufficiently small and apply the Vitali covering lemma to
\[
\{B_{r/5}(x):x\in\partial\Omega\}.
\]
There exist $x_1,\ldots,x_{N_r}\in\partial\Omega$ such that the balls
$B_{r/5}(x_k)$ are pairwise disjoint and
\[
\partial\Omega\subset\bigcup_{k=1}^{N_r}B_r(x_k).
\]
The enlarged balls $B_{2r}(x_k)$ have overlap bounded by a constant
depending only on $n$. By \eqref{admissible},
\[
\left|B_{2r}(x_k)\cap\{\widetilde m=0\}\right|
\geq |B_{2r}(x_k)\cap\Omega^c|
\geq\theta|B_{2r}(x_k)|.
\]
The Poincar\'e inequality therefore gives
\[
\int_{B_{2r}(x_k)}|\widetilde m|^2\,dx
\leq Cr^2\int_{B_{2r}(x_k)}|\nabla\widetilde m|^2\,dx.
\]

If $x\in\Omega^r$, choose $y\in\partial\Omega$ with $|x-y|<r$ and then
$x_k$ with $|y-x_k|<r$. Thus $x\in B_{2r}(x_k)$. Summing the preceding
inequalities and using bounded overlap yields
\[
\int_{\Omega^r}|m|^2\,dx
\leq Cr^2\int_{\Omega^{2r}}|\nabla m|^2\,dx.
\]
Since $|\nabla m|^2\in L^1(\Omega)$ and
$\Omega^{2r}\downarrow\varnothing$ as $r\downarrow0$,
\eqref{eq:boundary-layer-vanishing} follows.

Now let $z\in H^1(\mathbb R^n)$ vanish a.e. in $\Omega^c$. The same covering
argument proves \eqref{eq:boundary-poincare} with $z|_\Omega$ in place of
$m$. Choose $\xi_r\in W^{1,\infty}(\Omega)$ such that
\[
\xi_r=0\ \text{in }\Omega^r,\qquad
\xi_r=1\ \text{in }\Omega\setminus\Omega^{2r},\qquad
|\nabla\xi_r|\leq\frac Cr.
\]
Then $\xi_rz|_\Omega$ has compact support in $\Omega$, while
\[
\|(1-\xi_r)z\|_{H^1(\Omega)}^2
\leq
C\int_{\Omega^{2r}}\bigl(|z|^2+|\nabla z|^2\bigr)\,dx
+\frac C{r^2}\int_{\Omega^{2r}}|z|^2\,dx
\longrightarrow0.
\]
Hence $\xi_rz\to z$ strongly in $H^1(\Omega)$, and therefore
$z|_\Omega\in H_0^1(\Omega)$.
\end{proof}

\begin{rem}\label{rem:admissible-only}
The preceding argument uses only the exterior measure-density condition
\eqref{admissible}; it requires neither Lipschitz regularity nor a Sobolev
extension property of $\Omega$.
\end{rem}

\begin{lm}[Boundary truncation]\label{lm:boundary-truncation}
Let $\Omega$ be a bounded $\theta$-admissible domain and
$h\in H^1(\Omega)$. Suppose that there exist an open neighborhood $U$ of
$\partial\Omega$ and $\widetilde h\in H^1(U)$ such that
\[
\widetilde h=h\quad\text{a.e. in }U\cap\Omega,
\qquad
\widetilde h\leq0\quad\text{a.e. in }U\cap\Omega^c.
\]
Then $h^+\in H_0^1(\Omega)$.
\end{lm}

\begin{proof}
Choose $\chi\in C_0^\infty(U)$ such that $\chi=1$ in a neighborhood of
$\partial\Omega$. Since $\widetilde h^+\in H^1(U)$ and
$\widetilde h^+=0$ a.e. in $U\cap\Omega^c$, the zero extension of
$\chi\widetilde h^+$ from $U$ to $\mathbb R^n$ belongs to
$H^1(\mathbb R^n)$ and vanishes a.e. outside $\Omega$. Moreover,
$(1-\chi)h^+$ has compact support in $\Omega$, so its zero extension also
belongs to $H^1(\mathbb R^n)$. The sum of these two extensions is exactly
the zero extension of $h^+$ from $\Omega$ to $\mathbb R^n$. The last
assertion of Lemma \ref{lm:boundary-poincare} gives
$h^+\in H_0^1(\Omega)$.
\end{proof}

\begin{lm}[Equation-dependent truncation]\label{lm:equation-truncation}
Let $\Omega$ be a bounded $\theta$-admissible domain, and let
$w\in H^1_{\mathrm{loc}}(\Omega)$ be a weak subsolution of
\begin{equation}\label{eq:homogeneous-subsolution}
-D_i\!\left(a_{ij}(x)D_jw\right)\leq0
\quad\text{in }\Omega;
\end{equation}
that is,
\[
\int_\Omega a_{ij}(x)D_jwD_i\varphi\,dx\leq0
\]
for every nonnegative $\varphi\in C_0^\infty(\Omega)$. Suppose that there
exists an open neighborhood $V$ of $\partial\Omega$ and
$m\in H_0^1(\Omega)$ such that
\begin{equation}\label{eq:weak-boundary-control}
w\leq m\quad\text{a.e. in }V\cap\Omega.
\end{equation}
Then $w^+\in H_0^1(\Omega)$.
\end{lm}

\begin{proof}
Choose \(\rho_0>0\) such that
\[
\Omega^{4\rho_0}\subset V,
\qquad 2\rho_0<r_0,
\]
where \(r_0\) is given by Lemma \ref{lm:boundary-poincare}. For \(0<r<\rho_0\) and \(0<\delta<r/4\), define
\[
\eta_{\delta,r}(x)=
\begin{cases}
0,
& d(x,\partial\Omega)\leq\delta,\\[2mm]
\dfrac{d(x,\partial\Omega)-\delta}{\delta},
& \delta<d(x,\partial\Omega)<2\delta,\\[3mm]
1,
& 2\delta\leq d(x,\partial\Omega)\leq r,\\[2mm]
\dfrac{2r-d(x,\partial\Omega)}{r},
& r<d(x,\partial\Omega)<2r,\\[3mm]
0,
& d(x,\partial\Omega)\geq2r.
\end{cases}
\]
Since $\eta_{\delta,r}^2w^+$ has compact support in $\Omega$, we have
\[
\begin{aligned}
0&\geq\int_\Omega A\nabla w\cdot\nabla(\eta_{\delta,r}^2w^+)\,dx\\
&=\int_\Omega\eta_{\delta,r}^2 A\nabla w^+\cdot\nabla w^+\,dx
+2\int_\Omega\eta_{\delta,r}w^+ A\nabla w^+\cdot\nabla\eta_{\delta,r}\,dx\\
&\geq\frac1{2\Lambda}\int_\Omega\eta_{\delta,r}^2|\nabla w^+|^2\,dx
-C\int_\Omega|\nabla\eta_{\delta,r}|^2|w^+|^2\,dx.
\end{aligned}
\]
Thus

\[
\begin{aligned}
\int_{\{x\in\Omega:\,2\delta<d(x,\partial\Omega)<r\}}
|\nabla w^+|^2\,dx
&\leq
\frac{C}{\delta^2}\int_{\Omega^{2\delta}}|w^+|^2\,dx+
\frac{C}{r^2}\int_{\Omega^{2r}}|w^+|^2\,dx.
\end{aligned}
\]
Since \(w^+\leq m^+\) in \(\Omega^{2r}\), Lemma \ref{lm:boundary-poincare} yields
\[
\frac{1}{\delta^2}
\int_{\Omega^{2\delta}}|w^+|^2\,dx
\leq
\frac{1}{\delta^2}
\int_{\Omega^{2\delta}}|m|^2\,dx
\longrightarrow0
\qquad\text{as }\delta\downarrow0.
\]
Letting \(\delta\downarrow0\), we obtain
\[
\int_{\Omega^r}|\nabla w^+|^2\,dx
\leq
\frac{C}{r^2}\int_{\Omega^{2r}}|m|^2\,dx
\leq
C\int_{\Omega^{4r}}|\nabla m|^2\,dx.
\]
Together with \(w\in H^1_{\mathrm{loc}}(\Omega)\) and \(w^+\leq m^+\) near \(\partial\Omega\), this proves that \(w^+\in H^1(\Omega)\).

Now define
\[
\xi_r(x)=
\begin{cases}
0,
& d(x,\partial\Omega)\leq r,\\[2mm]
\dfrac{d(x,\partial\Omega)-r}{r},
& r<d(x,\partial\Omega)<2r,\\[3mm]
1,
& d(x,\partial\Omega)\geq2r.
\end{cases}
\]
Since \(\xi_rw^+\in H^1(\Omega)\) has compact support in \(\Omega\), it belongs to \(H_0^1(\Omega)\). Moreover,
\[
\begin{aligned}
\|\xi_rw^+-w^+\|_{H^1(\Omega)}^2
&\leq
C\int_{\Omega^{2r}}
\bigl(|w^+|^2+|\nabla w^+|^2\bigr)\,dx+
\frac{C}{r^2}\int_{\Omega^{2r}}|m|^2\,dx
\longrightarrow0
\qquad\text{as }r\downarrow0.
\end{aligned}
\]
Here the first term tends to zero because \(w^+\in H^1(\Omega)\), and the second tends to zero by Lemma \ref{lm:boundary-poincare}. Since \(H_0^1(\Omega)\) is closed in \(H^1(\Omega)\), we conclude that
\[
w^+\in H_0^1(\Omega).
\]
\end{proof}

\begin{rem}\label{rem:inhomogeneous-truncation}
Lemma \ref{lm:equation-truncation} remains valid if
\eqref{eq:homogeneous-subsolution} is replaced by
\[
-D_i\!\left(a_{ij}(x)D_jw\right)\leq f
\quad\text{in }\Omega,\qquad f\in L^2(\Omega).
\]
Moreover, after extension by zero outside $\Omega$, $w^+$ satisfies
\begin{equation}\label{eq:inhomogeneous-zero-extension}
\int_{\mathbb R^n}a_{ij}(x)D_jw^+D_i\varphi\,dx
\leq\int_\Omega f^+\varphi\,dx
\qquad\bigl(0\leq\varphi\in C_0^\infty(\mathbb R^n)\bigr).
\end{equation}
The test functions are not required to vanish on $\partial\Omega$.
The same test $\eta_{\delta,r}^2w^+$, followed by $\delta\downarrow0$, gives
\[
\begin{aligned}
\int_{\Omega^r}|\nabla w^+|^2\,dx
&\leq\frac C{r^2}\int_{\Omega^{2r}}|m|^2\,dx
+C\int_{\Omega^{2r}}|f|m^+\,dx\\
&\leq C\int_{\Omega^{4r}}|\nabla m|^2\,dx
+C\|f\|_{L^2(\Omega^{2r})}\|m\|_{L^2(\Omega^{2r})}
\longrightarrow0\quad(r\downarrow0).
\end{aligned}
\]
As in Lemma \ref{lm:equation-truncation}, $w^+\in H_0^1(\Omega)$.
For $0\leq\varphi\in C_0^\infty(\mathbb R^n)$ and $\varepsilon>0$,
\[
0\leq\xi_r\varphi\min\{w^+/\varepsilon,1\}\in H_0^1(\Omega),
\]
and the subsolution inequality gives
\[
\begin{aligned}
&\int_\Omega\xi_r\min\{w^+/\varepsilon,1\}
A\nabla w^+\cdot\nabla\varphi\,dx\\
&\quad+\int_\Omega\varphi\min\{w^+/\varepsilon,1\}
A\nabla w^+\cdot\nabla\xi_r\,dx\\
&\quad+\underbrace{\frac1\varepsilon\int_{\Omega\cap\{0<w<\varepsilon\}}
\xi_r\varphi A\nabla w^+\cdot\nabla w^+\,dx}_{\geq0}
\leq\int_\Omega f\xi_r\varphi\min\{w^+/\varepsilon,1\}\,dx.
\end{aligned}
\]
By Lemma \ref{lm:boundary-poincare},
\[
\begin{aligned}
&\left|\int_\Omega\varphi\min\{w^+/\varepsilon,1\}
A\nabla w^+\cdot\nabla\xi_r\,dx\right|\\
&\quad\leq\frac{C\|\varphi\|_{L^\infty}}{\varepsilon}
\left(\int_{\Omega^{2r}}|\nabla w^+|^2\,dx\right)^{1/2}
\left(\frac1{r^2}\int_{\Omega^{2r}}|m|^2\,dx\right)^{1/2}
\longrightarrow0\quad(r\downarrow0).
\end{aligned}
\]
Consequently,
\[
\begin{aligned}
\int_{\mathbb R^n}A\nabla w^+\cdot\nabla\varphi\,dx
&=\lim_{\varepsilon\downarrow0}\lim_{r\downarrow0}
\int_\Omega\xi_r\min\{w^+/\varepsilon,1\}
A\nabla w^+\cdot\nabla\varphi\,dx\\
&\leq\lim_{\varepsilon\downarrow0}
\int_\Omega f\varphi\min\{w^+/\varepsilon,1\}\,dx\\
&=\int_{\Omega\cap\{w>0\}}f\varphi\,dx
\leq\int_\Omega f^+\varphi\,dx.
\end{aligned}
\]

\end{rem}

\begin{lm}[Weak maximum principle]\label{lm:weak-maximum}
Under the assumptions of Lemma \ref{lm:equation-truncation},
\[
w\leq0\quad\text{a.e. in }\Omega.
\]
\end{lm}

\begin{proof}
By Lemma \ref{lm:equation-truncation}, $w^+\in H_0^1(\Omega)$.
Choose $r>0$ such that $\Omega^{2r}\subset V$, and let
$\xi_r\in W^{1,\infty}(\Omega)$ satisfy
\[
\xi_r=0\ \text{in }\Omega^r,\qquad
\xi_r=1\ \text{in }\Omega\setminus\Omega^{2r},\qquad
|\nabla\xi_r|\leq\frac Cr.
\]
Testing \eqref{eq:homogeneous-subsolution} with $\xi_r^2w^+$,
as in Lemma \ref{lm:equation-truncation}, gives

\[
\int_\Omega\xi_r^2|\nabla w^+|^2\,dx
\leq
\frac C{r^2}\int_{\Omega^{2r}}|w^+|^2\,dx.
\]
Since $w^+\leq m^+$ in $\Omega^{2r}$, Lemma
\ref{lm:boundary-poincare} yields
\[
\frac1{r^2}\int_{\Omega^{2r}}|w^+|^2\,dx
\leq
C\int_{\Omega^{4r}}|\nabla m|^2\,dx
\longrightarrow0
\quad\text{as }r\downarrow0.
\]
Because $\xi_r\to1$ pointwise in $\Omega$, Fatou's lemma implies
\[
\int_\Omega|\nabla w^+|^2\,dx=0.
\]
Finally, $w^+\in H_0^1(\Omega)$, so the Poincar\'e inequality gives
$w^+=0$ a.e. in $\Omega$.
\end{proof}

\begin{proof}[Proof of Theorem \ref{thm:uniqueness}]
Let $u$ and $v$ be two weak solutions with the same $f$ and $g$, and set
\[
w:=u-v.
\]
Subtracting their interior weak formulations gives
\[
-D_i\!\left(a_{ij}(x)D_jw\right)=0
\quad\text{in }\Omega
\]
in the weak sense.

Fix $\varepsilon>0$. Apply Proposition
\ref{prop:boundary-approximation} to $u$ and $v$, using the same
$U_\varepsilon$ and $g_\varepsilon^\pm$. In
$U_\varepsilon\cap\Omega$,
\[
u_\varepsilon^-\leq u\leq u_\varepsilon^+,
\qquad
v_\varepsilon^-\leq v\leq v_\varepsilon^+.
\]
Define
\[
W_\varepsilon^+:=u_\varepsilon^+-v_\varepsilon^-
\in H^1(\Omega).
\]
Then
\[
w\leq W_\varepsilon^+
\quad\text{a.e. in }U_\varepsilon\cap\Omega.
\]
Its extension
\[
\widetilde W_\varepsilon^+
:=\widetilde u_\varepsilon^+-\widetilde v_\varepsilon^-
\in H^1(U_\varepsilon)
\]
satisfies, by \eqref{eq:exterior-gap},
\[
\widetilde W_\varepsilon^+=\varepsilon
\quad\text{a.e. in }U_\varepsilon\cap\Omega^c.
\]
Applying Lemma \ref{lm:boundary-truncation} to
$W_\varepsilon^+-\varepsilon$ gives
\[
m_\varepsilon^+
:=\left(W_\varepsilon^+-\varepsilon\right)^+
\in H_0^1(\Omega).
\]
Moreover,
\[
w-\varepsilon
\leq W_\varepsilon^+-\varepsilon
\leq m_\varepsilon^+
\quad\text{a.e. in }U_\varepsilon\cap\Omega.
\]
Apply Lemma \ref{lm:weak-maximum} to $w-\varepsilon$, with
$m=m_\varepsilon^+$, to obtain
\[
w\leq\varepsilon\quad\text{a.e. in }\Omega.
\]
Interchanging $u$ and $v$ gives
\[
w\geq-\varepsilon\quad\text{a.e. in }\Omega.
\]
Letting $\varepsilon\downarrow0$ gives $u=v$ a.e. in $\Omega$.
\end{proof}

\section{Boundary pointwise \texorpdfstring{$C^\alpha$}{C-alpha} regularity}\label{sec:boundary-regularity}
We retain only the exterior measure-density condition
\eqref{admissible} to get  boundary $C^{\alpha}$ regularity compared with\cite{GWZ26}. Moreover, coefficients here are merely
measurable, bounded, and uniformly elliptic continuity wihout assuming VMO,
small-BMO, or
pointwise closeness-to-a-constant assumption on $A$. We restrict the available exponent to
$0<\alpha<\beta(n,\Lambda,\theta)$ as the price. In Section
\ref{sec:reifenberg-regularity}, Reifenberg flatness and small mean
oscillation of the coefficients allow any prescribed exponent
$0<\alpha<1$.
\subsection{Boundary H\"older regularity on admissible domains via harmonic replacements}

This part proves Theorem
\ref{thm:boundary-holder-harmonic-replacement}. It is independent of the
Reifenberg flatness and small-BMO assumptions.

We emphasize one point before the proof. In the classical Dirichlet theory,
one can take the $A$-harmonic replacement of $u$ in a subdomain $D$ by
requiring $h-u\in H_0^1(D)$. That construction cannot be used directly
here: the solution is known only to belong to $H^1_{\mathrm{loc}}(\Omega)$,
it need not belong to $H^1(D)$ up to the rough boundary, and no Sobolev trace
of $u$ on $\partial\Omega\cap\partial D$ is available. Thus the condition
$h-u\in H_0^1(D)$ may not even be meaningful. Definition \ref{uni} gives
locally $H^1$ zero extensions of the upper and lower boundary
truncations, which belong to $H^1$ on strictly smaller balls. We therefore form the classical $A$-harmonic replacements of
the explicit functions $(u-\varepsilon)^+$ and
$(u+\varepsilon)^-$, and use their one-sided subsolution inequalities to
compare them with those replacements.

\begin{proof}[Proof of Theorem
\ref{thm:boundary-holder-harmonic-replacement}]
Subtracting the constant $g(0)$, we may assume that $g(0)=0$.  Set
\[
 \Phi(r):=\left(
  \frac1{|B_r|}\int_{\Omega_r}|u|^2\,dx
 \right)^{1/2},
 \qquad
 K:=[g]_{C^\alpha(0)}+[f]_{\mathcal L^{2,-2+\alpha}(0)}.
\]

We first prove a comparison estimate entirely at unit scale.  Let
\[
 \varepsilon:=\sup_{\partial_w\Omega_1}|g|.
\]
Then
\begin{equation}
 0\le\varepsilon\le[g]_{C^\alpha(0)}.
 \label{eq:hm-epsilon}
\end{equation}
We use $(u-\varepsilon)^+$ and $(u+\varepsilon)^-$ to construct the
$A$-harmonic replacements. Indeed,
\[
\begin{aligned}
(u-\varepsilon)^+
&=\left((u-\sup_{\partial_w\Omega_1}g)^+
-\left(\varepsilon-\sup_{\partial_w\Omega_1}g\right)\right)^+,\\
(u+\varepsilon)^-
&=\left((u-\inf_{\partial_w\Omega_1}g)^-
-\left(\varepsilon+\inf_{\partial_w\Omega_1}g\right)\right)^+.
\end{aligned}
\]
Both constants subtracted on the right are nonnegative. Lemma
\ref{lm:further-truncation} therefore shows that the zero extensions of
$(u-\varepsilon)^+$ and $(u+\varepsilon)^-$ belong to $H^1_{\mathrm{loc}}(B_1)$, vanish
almost everywhere in $B_1\setminus\Omega$, and satisfy
\begin{align}
 \int_{B_1}a_{ij}(x)D_j[(u-\varepsilon)^+]D_i\varphi\,dx
 &\le\int_{B_1}f^+\varphi\,dx,
 \label{eq:hm-subsolution-plus}\\
 \int_{B_1}a_{ij}(x)D_j[(u+\varepsilon)^-]D_i\varphi\,dx
 &\le\int_{B_1}f^-\varphi\,dx
 \label{eq:hm-subsolution-minus}
\end{align}
for every nonnegative function $\varphi\in C_0^{\infty}(B_1)$. 

Choose $\eta\in C_c^\infty(B_1)$ such that $0\le\eta\le1$ and
$\eta=1$ on $B_{1/2}$.  Applying the boundary energy argument of Lemma
\ref{lm1} to these further truncations gives the following boundary Caccioppoli
estimates
\begin{equation}
\begin{aligned}
 \int_{\Omega_{1/2}}|\nabla (u-\varepsilon)^+|^2\,dx
 &\le C\int_{\Omega_1}\bigl(|u|^2+|f|^2\bigr)\,dx,\\
 \int_{\Omega_{1/2}}|\nabla (u+\varepsilon)^-|^2\,dx
 &\le C\int_{\Omega_1}\bigl(|u|^2+|f|^2\bigr)\,dx.
\end{aligned}
 \label{eq:hm-caccioppoli}
\end{equation}
Here we used
$0\le (u-\varepsilon)^+,(u+\varepsilon)^-\le|u|$.

Set $D:=\Omega_{1/2}$.  Let $h^\pm\in H^1(D)$ be the classical
$A$-harmonic replacements of $(u-\varepsilon)^+$ and
$(u+\varepsilon)^-$, respectively, namely
\begin{equation}
\begin{gathered}
 h^+-(u-\varepsilon)^+\in H_0^1(D),\qquad
 \int_D a_{ij}(x)D_jh^+D_i\varphi\,dx=0
 \quad\forall\varphi\in H_0^1(D),\\
 h^--(u+\varepsilon)^-\in H_0^1(D),\qquad
 \int_D a_{ij}(x)D_jh^-D_i\varphi\,dx=0
 \quad\forall\varphi\in H_0^1(D).
\end{gathered}
 \label{eq:hm-replacement}
\end{equation}
This makes sense because $(u-\varepsilon)^+$ and $(u+\varepsilon)^-$ belong to $H^1(D)$. The weak maximum principle gives $h^\pm\ge0$.

The zero extensions of $h^\pm$ belong to $H^1(B_{1/2})$ and vanish
almost everywhere in $B_{1/2}\setminus\Omega$.  To see this, write
\[
\begin{aligned}
 h^+&=(u-\varepsilon)^++\bigl(h^+-(u-\varepsilon)^+\bigr),\\
 h^-&=(u+\varepsilon)^-+\bigl(h^--(u+\varepsilon)^-\bigr),
\end{aligned}
\]
and extend the second term by zero.  This follows directly from
\eqref{eq:hm-replacement}.

Taking $h^+-(u-\varepsilon)^+\in H_0^1(D)$ in
\eqref{eq:hm-replacement}, we have
\[
\begin{aligned}
\Lambda^{-1}\|\nabla(h^+-(u-\varepsilon)^+)\|_{L^2(D)}^2
&\leq-\int_D A\nabla(u-\varepsilon)^+\cdot
\nabla(h^+-(u-\varepsilon)^+)\,dx\\
&\leq\Lambda\|\nabla(u-\varepsilon)^+\|_{L^2(D)}
\|\nabla(h^+-(u-\varepsilon)^+)\|_{L^2(D)}.
\end{aligned}
\]
The same calculation for $h^--(u+\varepsilon)^-$, Poincar\'e's
inequality, and \eqref{eq:hm-caccioppoli} give

\begin{equation}
 \|h^\pm\|_{L^2(D)}
 \le C\left(
  \|u\|_{L^2(\Omega_1)}+\|f\|_{L^2(\Omega_1)}
 \right).
 \label{eq:hm-replacement-L2}
\end{equation}

The zero-extension form of the boundary De Giorgi--Moser estimate now
yields, for some
$\beta=\beta(n,\Lambda,\theta)\in(0,1)$ and every
$0<\rho\le1/8$,
\begin{equation}
 \left(
  \frac1{|B_\rho|}\int_{\Omega_\rho}|h^\pm|^2\,dx
 \right)^{1/2}
 \le C\rho^\beta
 \left(
  \frac1{|B_{1/2}|}\int_D|h^\pm|^2\,dx
 \right)^{1/2}.
 \label{eq:hm-moser}
\end{equation}

In \eqref{eq:hm-subsolution-plus}--\eqref{eq:hm-subsolution-minus},
use the zero extensions of
\[
\bigl((u-\varepsilon)^+-h^+\bigr)^+,
\quad\bigl((u+\varepsilon)^--h^-\bigr)^+\in H_0^1(D).
\]
Subtracting \eqref{eq:hm-replacement} gives
\[
\begin{aligned}
&\Lambda^{-1}\int_D
\bigl|\nabla\bigl((u-\varepsilon)^+-h^+\bigr)^+\bigr|^2\,dx\\
&\quad\leq\int_D A\nabla\bigl((u-\varepsilon)^+-h^+\bigr)\cdot
\nabla\bigl((u-\varepsilon)^+-h^+\bigr)^+\,dx
\leq\int_D f^+\bigl((u-\varepsilon)^+-h^+\bigr)^+\,dx,\\
&\Lambda^{-1}\int_D
\bigl|\nabla\bigl((u+\varepsilon)^--h^-\bigr)^+\bigr|^2\,dx\\
&\quad\leq\int_D A\nabla\bigl((u+\varepsilon)^--h^-\bigr)\cdot
\nabla\bigl((u+\varepsilon)^--h^-\bigr)^+\,dx
\leq\int_D f^-\bigl((u+\varepsilon)^--h^-\bigr)^+\,dx.
\end{aligned}
\]
Poincar\'e's inequality therefore gives
\begin{equation}
\begin{aligned}
 \bigl\|\nabla\bigl((u-\varepsilon)^+-h^+\bigr)^+\bigr\|_{L^2(D)}
 +\bigl\|\bigl((u-\varepsilon)^+-h^+\bigr)^+\bigr\|_{L^2(D)}
 &\le C\|f^+\|_{L^2(D)},\\
 \bigl\|\nabla\bigl((u+\varepsilon)^--h^-\bigr)^+\bigr\|_{L^2(D)}
 +\bigl\|\bigl((u+\varepsilon)^--h^-\bigr)^+\bigr\|_{L^2(D)}
 &\le C\|f^-\|_{L^2(D)}.
\end{aligned}
 \label{eq:hm-one-sided}
\end{equation}
Notice that the above letter $C$ depends on $R$ when we calculate on $\Omega\cap B_{R}$.  In particular,
\begin{equation}
\begin{aligned}
 (u-\varepsilon)^+
 &\le h^++\bigl((u-\varepsilon)^+-h^+\bigr)^+,\\
 (u+\varepsilon)^-
 &\le h^-+\bigl((u+\varepsilon)^--h^-\bigr)^+.
\end{aligned}
 \label{eq:hm-upper}
\end{equation}

For every real number $u$,
\[
 |u|\le\varepsilon+(u-\varepsilon)^++(u+\varepsilon)^-.
\]
Using \eqref{eq:hm-upper}, followed by
\eqref{eq:hm-replacement-L2}--\eqref{eq:hm-one-sided}, gives the unit-scale
estimate
\begin{equation}
 \Phi(\rho)
 \le C\rho^\beta\Phi(1)+C\varepsilon
 +C\bigl(\rho^\beta+\rho^{-n/2}\bigr)
 \left(
  \frac1{|B_1|}\int_{\Omega_1}|f|^2\,dx
 \right)^{1/2}
 \label{eq:hm-unit-scale}
\end{equation}
for $0<\rho\le1/8$.

To rescale \eqref{eq:hm-unit-scale}, fix $0<R\leq1$ and set
\[
\widetilde\Omega:=R^{-1}\Omega,\qquad
u_R(x):=u(Rx),\qquad A_R(x):=A(Rx),
\]
\[
f_R(x):=R^2f(Rx),\qquad g_R(x):=g(Rx).
\]
Then $u_R$ satisfies the same weak formulation on $\widetilde\Omega$,
with coefficient matrix $A_R$, forcing $f_R$, and boundary data $g_R$.
The ellipticity and exterior density constants are unchanged.
The change of variables $y=Rx$ gives, for $0<t\leq1$,
\[
\left(\frac1{|B_t|}
\int_{\widetilde\Omega\cap B_t}|u_R|^2\,dx\right)^{1/2}
=\left(\frac1{|B_{Rt}|}\int_{\Omega_{Rt}}|u|^2\,dy\right)^{1/2}
=\Phi(Rt),
\]
and
\[
\left(\frac1{|B_1|}
\int_{\widetilde\Omega\cap B_1}|f_R|^2\,dx\right)^{1/2}
=R^2\left(\frac1{|B_R|}\int_{\Omega_R}|f|^2\,dy\right)^{1/2}.
\]
Applying \eqref{eq:hm-unit-scale} to $u_R$ with $\rho=\tau$, where
$0<\tau\leq1/8$, therefore yields
\[
\begin{aligned}
\Phi(\tau R)
&\leq C\tau^\beta\Phi(R)
+C\sup_{\partial_w\Omega_R}|g|\\
&\quad+C\bigl(\tau^\beta+\tau^{-n/2}\bigr)R^2
\left(\frac1{|B_R|}\int_{\Omega_R}|f|^2\,dy\right)^{1/2}.
\end{aligned}
\]
Since $g(0)=0$, the definitions of the pointwise seminorms give
\[
\sup_{\partial_w\Omega_R}|g|
\leq [g]_{C^\alpha(0)}R^\alpha,\qquad
R^2\left(\frac1{|B_R|}\int_{\Omega_R}|f|^2\,dy\right)^{1/2}
\leq [f]_{\mathcal L^{2,-2+\alpha}(0)}R^\alpha.
\]
Recalling $K=[g]_{C^\alpha(0)}+[f]_{\mathcal L^{2,-2+\alpha}(0)}$
and including the fixed factors involving $\tau$ in $C_\tau$, we obtain
\begin{equation}
 \Phi(\tau R)
 \le C\tau^\beta\Phi(R)+C_\tau K R^\alpha.
 \label{eq:hm-scaled}
\end{equation}

Choose $\tau$ so small that
\[
 C\tau^\beta\le\frac12\tau^\alpha.
\]
By \eqref{eq:hm-scaled},
\[
\begin{aligned}
\tau^{-(k+1)\alpha}\Phi(\tau^{k+1})
&\leq\frac12\tau^{-k\alpha}\Phi(\tau^k)+C_\tau\tau^{-\alpha}K,\\
\tau^{-k\alpha}\Phi(\tau^k)
&\leq2^{-k}\Phi(1)+C_\tau\tau^{-\alpha}K\sum_{j=0}^{k-1}2^{-j}
\leq\Phi(1)+2C_\tau\tau^{-\alpha}K.
\end{aligned}
\]
For arbitrary $0<r\le1/8$, choose $k$ so that
$\tau^{k+1}<r\le\tau^k$.  Since
$\Omega_r\subset\Omega_{\tau^k}$,
\[
 \Phi(r)
 \le\left(\frac{\tau^k}{r}\right)^{n/2}\Phi(\tau^k)
 \le Cr^\alpha\bigl(\Phi(1)+K\bigr).
\]
Restoring $g(0)$ proves the theorem.
\end{proof}

\section{Regularity on Reifenberg domains}\label{sec:reifenberg-regularity}
We first prove boundary pointwise $C^\alpha$ regularity for any prescribed
$0<\alpha<1$ by the compactness method, under Reifenberg flatness and
the small-BMO assumption on $A$. We then prove global $H^1$ regularity
by constructing a continuous energy lifting of the boundary data and
applying uniqueness. The second argument does not require small mean
oscillation of the coefficients.
\subsection{Boundary pointwise \texorpdfstring{$C^\alpha$}{C-alpha} regularity}
We show the following approximation lemma by the compactness method.
\begin{lm}[Compactness lemma]\label{lm4}
For any $\varepsilon>0$, there exists a small $\delta=\delta(\varepsilon,n,\Lambda)>0$ such that, if $A$ satisfies \eqref{ellipticity}, then for any weak solution of $(\ref{1})$ in $\Omega_{1}$ with $$\displaystyle\frac{1}{|B_{1}|}\displaystyle\int_{\Omega_{1}}u^{2}dx\leq 1, ~\frac 1{|B_1|}\int_{\Omega_{1}}f^{2}dx\leq \delta^2, ~\|g\|_{L^{\infty}(\partial_w\Omega_1)}\leq \delta,~\frac{1}{|B_1|}\int_{\Omega_{1}}
\big|A-(A)_{\Omega_1}\big|^2\,dx \leq \delta^{2},$$
and $$B_1^{+}\subset\Omega_{1}\subset B_1\cap \{x_n>-\delta\},$$
there exist a constant matrix $\overline A=(\overline a_{ij})$ satisfying \eqref{ellipticity} and a corresponding weak solution $h$ of
\[
-D_i\!\left(\overline a_{ij}D_jh\right)=0
\quad\text{in }B^+_{\frac12},
\qquad h=0\quad\text{on }T_{\frac12},
\]
such that
$$\int_{B_{\frac{1}{2}}^{+}}|u-h|^{2}dx\leq\varepsilon^{2}.
$$
\end{lm}
\begin{proof}
We prove it by contradiction. If not, there is an $\varepsilon_0 >0$ so that for each integer $k$, there is a weak solution of
\begin{eqnarray*}
\left\{
\begin{array}{rcll}
-D_i\!\left(a^k_{ij}(x)D_j u_k\right)&=&f_k\qquad&\text{in}~~\Omega_{k},\\
u_k&=&g_k\qquad&\text{on}~~\partial\Omega_{k}, \\
\end{array}
\right.
\end{eqnarray*}
where each $A_k=(a^k_{ij})$ satisfies \eqref{ellipticity} and
\[
\begin{gathered}
\frac{1}{|B_{1}|}\int_{(\Omega_k)_1}u_{k}^{2}\,dx\leq 1,
\qquad
\frac{1}{|B_1|}\int_{(\Omega_k)_1}f_{k}^{2}\,dx\leq \frac{1}{k^2},
\qquad
\|g_k\|_{L^{\infty}(\partial_w(\Omega_k)_1)}\leq \frac{1}{k},\\
\frac{1}{|B_1|}\int_{(\Omega_{k})_1}
\big|A_{k}-(A_k)_{(\Omega_k)_1}\big|^2\,dx \leq \frac{1}{k^2},
\end{gathered}
\]
and
$$B_{1}^{+}\subset\Omega_{k}\subset B_1\cap \{x_n>-\frac1k\},$$
 such that for any constant matrix $\overline A=(\overline a_{ij})$ satisfying \eqref{ellipticity} and any solution $h$ of $-D_i(\overline a_{ij}D_jh)=0$ in $B^+_{\frac12}$ with zero data on $T_{\frac12}$, 
 the following inequality holds
\begin{equation}\label{c}\int_{B^{+}_{\frac{1}{2}}}|u_{k}-h|^{2}dx>\varepsilon_{0}^{2}.
\end{equation}

%From global $C^{\alpha}$ estimate, we have
%\begin{equation}\label{gl1}
%    (u_k-\sup_{\partial_w(\Omega_k)_1}g_k)^{+}\leq C[d^{\alpha}(x,\partial\Omega_k)+\|f_k\|_{L^{2}%(\Omega_k)}]
%\end{equation}
%\begin{equation}\label{gl2}
%    (u_k-\inf_{\partial_w(\Omega_k)_1}g_k)^{-}\leq C[d^{\alpha}(x,\partial\Omega_k)+\|f_k\|_{L^{2}%(\Omega_k)}]
%\end{equation}
From Lemma $\ref{lm1}$ we observe that
$$\|(u_k-\sup\limits_{\partial_w(\Omega_k)_1}g_k)^{+}\|_{H^1(B_{\frac{3}{4}})}\leq C,$$
$$\|(u_k-\inf\limits_{\partial_w(\Omega_k)_1}g_k)^{-}\|_{H^1(B_{\frac{3}{4}})}\leq C.$$

Set $\overline A_k=(A_k)_{(\Omega_k)_1}$. By ellipticity and boundedness, after
passing to a subsequence, $\overline A_k\to\overline A$ for a constant matrix
$\overline A$ satisfying \eqref{ellipticity}. The last smallness assumption
in the contradiction sequence gives
\[
\|A_k-\overline A_k\|_{L^2((\Omega_k)_1)}
\leq \frac{|B_1|^{1/2}}{k},
\]
and hence $A_k\to\overline A$ strongly in $L^2$ on every compact subset of
$B_1^+$.

Consequently, using the diagonal subsequence technique, it is easy to see that $(u_k-\sup\limits_{\partial_w(\Omega_k)_1}g_k)^{+}$ [$(u_k-\inf\limits_{\partial_w(\Omega_k)_1}g_k)^{-}$, resp.] has a subsequence, which we still denoted as $(u_k-\sup\limits_{\partial_w(\Omega_k)_1}g_k)^{+}$ [$(u_k-\inf\limits_{\partial_w(\Omega_k)_1}g_k)^{-}$, resp.] that the weak limit $v_1$[$v_2$, resp.] $\in H^{1}(B_{\frac{3}{4}})$ and strongly converges to $v_1$[$v_2$,  resp.] in $L^{2}(B_{\frac{3}{4}})$.  However, we can only  deduce that $u_k$ weakly converges to $u_{\infty}$ in $L^{2}(B_{\frac{3}{4}})$. Our aim is to prove that $u_k$ converges to $u_{\infty}$ strongly in $L^{2}(B_{\frac{3}{4}})$ and $u_{\infty}^{+}=v_1 $, $u_{\infty}^{-}=v_2 $ a.e.. Then we have $u_\infty=v_1-v_2$ a.e..

We write $u_k$ as $u_k^{+}-u_k^{-}$, and let $u_k=0$ in $((\Omega_k)^c)^o$. Notice that $(u_k-\sup\limits_{\partial_w(\Omega_k)_1}g_k)^{+}$ is a Cauchy sequence in $L^{2}(B_{\frac{3}{4}})$, so is $u_{k}^{+}$. In fact, we have
\begin{eqnarray*}\|u^{+}_{k+t}-u^{+}_{k}\|_{L^{2}(B_{\frac{3}{4}})}
&\leq& \|u^{+}_{k+t}-(u_{k+t}-\sup\limits_{\partial_w(\Omega_{k+t})_1}g_{k+t})^{+}\|_{L^{2}(B_{\frac{3}{4}})}\\
& & +\|(u_{k+t}-\sup\limits_{\partial_w(\Omega_{k+t})_1}g_{k+t})^{+}-(u_{k}-\sup\limits_{\partial_w(\Omega_k)_1}g_{k})^{+}\|_{L^{2}(B_{\frac{3}{4}})}
\\
& & +\|(u_{k}-\sup\limits_{\partial_w(\Omega_k)_1}g_{k})^{+}-u^{+}_{k}\|_{L^{2}(B_{\frac{3}{4}})} \\
&\rightarrow& 0\quad \text{as} \quad k,t\rightarrow \infty
\end{eqnarray*}
From above we can assume that $u_k^+$ converges to $u_{\infty}^+$ strongly in $L^{2}(B_{\frac{3}{4}})$. Since
$$|u^{+}_{\infty}-v_1|
\leq |u^{+}_{\infty}-u^{+}_{k}|+| u^+_{k}- (u_k-\sup\limits_{\partial_w(\Omega_k)_1}g_{k})^{+}|+| (u_k-\sup\limits_{\partial_w(\Omega_k)_1}g_{k})^{+}-v_1|,
$$
We have $u_{\infty}^{+}=v_1 $ a.e.. Similarly we have
$u_k^-$ converges to $u_{\infty}^-$ strongly in $L^{2}(B_{\frac{3}{4}})$ and $u_{\infty}^{-}=v_2 $ a.e. .

Writing $u_{\infty}$ as $v$, since $u_k=0$ outside $\Omega_k$ and
$\Omega_k\subset B_1\cap\{x_n>-1/k\}$, we also have
$v=0$ a.e. in $B^-_{\frac12}$.
The preceding local $H^1$ bounds also give, after passing to the same
subsequence, $\nabla u_k\rightharpoonup \nabla v$ weakly in $L^2$ on compact subsets of
$B^+_{\frac12}$.

For any $ x_0\in  B^{+}_{\frac{1}{2}}$, there exists $k_0$ and $t_0>0$, such that for any $k\geq k_0$, $x_0\in \Omega_k$ and $B_{t_0}(x_0)\subset \Omega_k$. Then for any $\eta \in C_{0}^{\infty}(B_{t_{0}}(x_0))$,
$$\int_{B_{t_{0}}(x_0)}f_k\eta\,dx
=\int_{B_{t_{0}}(x_0)}a^k_{ij}(x)D_ju_kD_i\eta\,dx.$$
Letting $k\rightarrow\infty$, using the weak convergence of $\nabla u_k$ and the
strong $L^2$ convergence of $A_k$, we obtain
\[
\int_{B_{t_0}(x_0)}\overline a_{ij}D_jvD_i\eta\,dx=0.
\]
Thus $v$ is a weak solution of $-D_i(\overline a_{ij}D_jv)=0$
in $B_{t_0}(x_0)$, and hence in $B^+_{\frac 12}$.
By the boundary trace argument in the proof of Lemma 4.5 in
\cite{BW04}, $v=0$ on $T_{\frac12}$. Finally, we get a contradiction to \eqref{c} by taking $h=v$, and $k$ large enough.
\end{proof}

Next we give a key lemma, which will be used repeatedly in the iteration.
\begin{lm}[Key Lemma]\label{key22} For any $0<\alpha<1$, there exists $0<\lambda<1,~\delta^{*}=\delta^*(\alpha,n,\Lambda)>0$ such that, if $A$ satisfies \eqref{ellipticity}, then for any weak solution of $(\ref{1})$
with
$$\displaystyle{\frac{1}{|B_1|}}\int_{\Omega_{1}}u^{2}dx\leq 1,~~
\frac 1{|B_1|}\int_{\Omega_{1}}f^{2}dx\leq (\delta^{*})^2,  ~~ \|g\|_{L^{\infty}(\partial_w\Omega_1)}\leq \delta^*,\frac{1}{|B_1|}\int_{\Omega_{1}}
\big|A-(A)_{\Omega_1}\big|^2\,dx \leq (\delta^{*})^{2}$$ 
and$$B_1^{+}\subset\Omega_{1}\subset B_1\cap \{x_n>-\delta^*\},$$
it follows that
$$\left(\frac{1}{|B_{\lambda}|}\int_{\Omega_{\lambda}}|u|^{2}dx\right)^{\frac{1}{2}}\leq\lambda^{\alpha}.
$$

\end{lm}
\begin{proof}
Let $h$ be the function given by Lemma \ref{lm4}, and let $H$ be its zero
extension from $B^+_{\frac12}$ to $B_{\frac12}$. Then
\begin{equation}\label{eq:key-compactness}
\int_{B^+_{\frac12}}|u-h|^2\,dx\leq\varepsilon^2
\end{equation}
for some $0<\varepsilon<1$ to be chosen below.

We first estimate $u-H$ on the part of $\Omega_{\frac12}$ below the flat
boundary. The zero extensions of
$\left(u-\sup_{\partial_w\Omega_1}g\right)^+$ and
$\left(u-\inf_{\partial_w\Omega_1}g\right)^-$ belong to $H^1_{\mathrm{loc}}(B_1)$.
Lemma \ref{lm1}, applied with a cutoff which is one on $B_{\frac34}$, and
the normalization assumptions give
\begin{equation}\label{eq:key-truncation-H1}
\left\|\left(u-\sup_{\partial_w\Omega_1}g\right)^+\right\|_{H^1(B_{\frac34})}
+\left\|\left(u-\inf_{\partial_w\Omega_1}g\right)^-\right\|_{H^1(B_{\frac34})}
\leq C.
\end{equation}
Let
\[
E_{\delta^*}:=\Omega_{\frac12}\setminus B^+_{\frac12}.
\]
Since
\[
E_{\delta^*}\subset
B_{\frac12}\cap\{-\delta^*<x_n\leq0\},
\qquad |E_{\delta^*}|\leq C\delta^*,
\]
and $\|g\|_{L^\infty(\partial_w\Omega_1)}\leq\delta^*$, we have
\[
|u|\leq
\left(u-\sup_{\partial_w\Omega_1}g\right)^+
+\left(u-\inf_{\partial_w\Omega_1}g\right)^-
+2\delta^*
\qquad\text{a.e. in }E_{\delta^*}.
\]
For $n>2$, take the Sobolev exponent
\[
q=2^*:=\frac{2n}{n-2}.
\]
By H\"older's inequality, the Sobolev inequality, and
\eqref{eq:key-truncation-H1},
\begin{align}
\int_{E_{\delta^*}}|u|^2\,dx
&\leq C|E_{\delta^*}|^{\frac2n}
\left(
\left\|\left(u-\sup_{\partial_w\Omega_1}g\right)^+\right\|_{L^{\frac{2n}{n-2}}(B_{\frac34})}^2
+\left\|\left(u-\inf_{\partial_w\Omega_1}g\right)^-\right\|_{L^{\frac{2n}{n-2}}(B_{\frac34})}^2
\right)
+C(\delta^*)^2|E_{\delta^*}| \notag\\
&\leq C(\delta^*)^{\frac2n}. \label{eq:key-thin-layer}
\end{align}
When $n=2$, repeat the estimates below with any fixed $q>2$, replacing
$2/n$ by $1-2/q$ and $1/n$ by $(1-2/q)/2$.
Because $H=0$ on $E_{\delta^*}$, \eqref{eq:key-compactness} and
\eqref{eq:key-thin-layer} yield
\begin{equation}\label{eq:key-u-minus-H}
\|u-H\|_{L^2(\Omega_{\frac12})}
\leq\left(\varepsilon^2+C(\delta^*)^{\frac2n}\right)^{1/2}.
\end{equation}

We next estimate $H$. From \eqref{eq:key-compactness} and the normalization
of $u$,
\[
\|h\|_{L^2(B^+_{\frac12})}
\leq\|u-h\|_{L^2(B^+_{\frac12})}
+\|u\|_{L^2(B^+_{\frac12})}\leq C.
\]
Since $\overline A$ is constant, its skew-symmetric part does not
contribute to $-D_i(\overline a_{ij}D_jh)$. The symmetric part satisfies
the same coercivity and boundedness assumptions. The local boundary
$W^{1,\infty}$ estimate for the constant-coefficient equation,
together with $h=0$ on $T_{\frac12}$, gives
\begin{equation}\label{eq:key-H-Lipschitz}
\|H\|_{W^{1,\infty}(B_{\frac14})}\leq C.
\end{equation}
In particular, $H(0)=0$, and hence, for $0<\lambda<\frac14$,
\begin{equation}\label{eq:key-H-L2}
\|H\|_{L^2(\Omega_\lambda)}
\leq C\lambda |B_\lambda|^{1/2}.
\end{equation}

Finally, using $u=(u-H)+H$, \eqref{eq:key-u-minus-H}, and
\eqref{eq:key-H-L2}, we obtain
\begin{align*}
\left(\frac1{|B_\lambda|}
\int_{\Omega_\lambda}|u|^2\,dx\right)^{1/2}
&\leq \frac{\|u-H\|_{L^2(\Omega_\lambda)}}{|B_\lambda|^{1/2}}
+\frac{\|H\|_{L^2(\Omega_\lambda)}}{|B_\lambda|^{1/2}}\\
&\leq
\frac{\left(\varepsilon^2+C(\delta^*)^{\frac2n}\right)^{1/2}}
{|B_\lambda|^{1/2}}+C\lambda.
\end{align*}
Choose $0<\lambda<\frac14$, depending only on $\alpha,n,\Lambda$,
so that $C\lambda\leq\frac12\lambda^\alpha$. Next choose $\varepsilon>0$
so that
\[
\frac{\varepsilon}{|B_\lambda|^{1/2}}
\leq\frac14\lambda^\alpha,
\]
and finally choose $\delta^*>0$ so that
\[
\delta^*\leq\delta(\varepsilon,n,\Lambda)
\quad\text{from Lemma \ref{lm4}},
\qquad
\frac{C^{1/2}(\delta^*)^{\frac1n}}{|B_\lambda|^{1/2}}
\leq\frac14\lambda^\alpha.
\]
The desired estimate follows.
\end{proof}

{\bf{Proof of Theorem $\ref{thm1}$:}}

We first discuss the normalization of the estimates. Without loss of
generality, we may assume $g(0)=0$. Otherwise, let $v(x)=u(x)-g(0)$.
Then the estimate of $v$ will translate to that for $u$.

Let $\delta^*$ be as in Lemma $\ref{key22}$. We can assume the
normalization conditions
\[
\begin{gathered}
\frac{1}{|B_1|}\int_{\Omega_1}u^2\,dx\leq1,
\qquad
[f]^2_{\mathcal L^{2,-2+\alpha}(0)}\leq(\delta^*)^2,
\qquad
[g]_{C^\alpha(0)}\leq\delta^*,\\
\frac1{|B_1|}\int_{\Omega_1}
\big|A-(A)_{\Omega_1}\big|^2\,dx\leq(\delta^*)^2.
\end{gathered}
\]
Otherwise, we set
$$
v(x)=\frac{u(x)}{
\left(\frac{1}{|B_1|}\int_{\Omega_1}u^2\,dx\right)^{1/2}
+\frac{1}{\delta^*}\left(
[f]_{\mathcal L^{2,-2+\alpha}(0)}+[g]_{C^\alpha(0)}\right)}
\triangleq\frac{u(x)}{N},\qquad x\in B_1\cap\Omega.
$$
Then $v$ is a weak solution of
\begin{equation}\label{4-2}
\left\{
\begin{array}{rcll}
-D_i\!\left(a_{ij}(x)D_j v\right)&=&{f}_{N}\qquad&\text{in}~~\Omega,\\
v&=&g_{N}\qquad&\text{on}~~\partial\Omega,\\
\end{array}
\right.
\end{equation}
where $f_N=f/N$ and $g_N=g/N$. Moreover,
\[
\begin{gathered}
\left(\frac{1}{|B_1|}\int_{\Omega_1}v^2\,dx\right)^{1/2}\leq1,
\qquad
[f_N]^2_{\mathcal L^{2,-2+\alpha}(0)}\leq(\delta^*)^2,
\qquad
[g_N]_{C^\alpha(0)}\leq\delta^*,\\
\frac1{|B_1|}\int_{\Omega_1}
\big|A-(A)_{\Omega_1}\big|^2\,dx\leq(\delta^*)^2.
\end{gathered}
\]

Now, we prove the following claim inductively:
\begin{equation}\label{ite3}
\left(\frac{1}{|B_{\lambda^{k}}|}\int_{\Omega_{\lambda^{k}}}u^{2}dx\right)^{\frac{1}{2}}\leq \lambda^{k\alpha},\quad \forall k\geq0.
\end{equation}
When $k=0$, (\ref{ite3}) is the condition on $u$. By normalized assumption, we know that
$$
\frac 1{|B_1|}\int_{\Omega_1}|f(x)|^2dx\leq [f]^2_{\mathcal{L}^{2,-2+\alpha}(0)}\leq (\delta^{*})^2,
$$
and, since $g(0)=0$,
$$
\|g\|_{L^{\infty}(\partial_w\Omega_1)}\leq [g]_{C^{\alpha}(0)}\leq\delta^{*}.
$$
Moreover, since $A$ is $(\delta,1)$-vanishing and $\delta<\delta^*$,
\[
\frac1{|B_1|}\int_{\Omega_1}|A-(A)_{\Omega_1}|^2\,dx
\leq\delta^2< (\delta^*)^2.
\]
Hence $(\ref{ite3})$  holds for $k=1$  from Lemma $\ref{key22}$.
Assume that the conclusion is true for some $k\geq1$. Let
$$w(x)=\frac{u(\lambda^{k}x)}{\lambda^{k\alpha}},\quad
\widetilde{\Omega}:=\{x\in\mathbb R^n:\lambda^k x\in\Omega\}
=\lambda^{-k}\Omega,\quad
\widetilde a_{ij}(x):=a_{ij}(\lambda^k x),\quad
\widetilde A(x):=(\widetilde a_{ij}(x)).
$$
It is easy to check that for any $k\geq 0$, $\widetilde{\Omega}$ is
$(\delta,1)$-Reifenberg flat, and $w$ is a weak solution to
\begin{eqnarray*}
\left\{
\begin{array}{rcll}
-D_i\!\left(\widetilde a_{ij}(x)D_jw\right)&=&\widetilde{f}\qquad&\text{in}~~\widetilde{\Omega},\\
w&=&\widetilde{g}\qquad&\text{on}~~\partial\widetilde{\Omega}, \\
\end{array}
\right.
\end{eqnarray*}
where $\widetilde{f}(x)
=\displaystyle\frac{f(\lambda^{k}x)}{\lambda^{k(\alpha-2)}}$, $\widetilde{g}(x)
=\displaystyle\frac{g(\lambda^{k}x)}{\lambda^{k\alpha}}$.
Thus by using the inductive assumption, we obtain that
$$\frac{1}{|B_{1}|}\displaystyle\int_{\widetilde{\Omega}_{1}}w^{2}dx=
\frac{1}{|B_{1}|}\displaystyle\int_{\widetilde{\Omega}_{1}}\frac{|u(\lambda^{k}x)|^{2}}{\lambda^{2k\alpha}}dx
=\frac{\displaystyle\frac{1}{|B_{\lambda^{k}}|}\displaystyle\int_{\Omega_{\lambda^{k}}}|u|^{2}dx}{\lambda^{2k\alpha}}
\leq1,
$$
$$\frac 1{|B_1|}\int_{\widetilde{\Omega}_{1}}\widetilde{f}^{2}dx=\frac 1{|B_1|}
\int_{\widetilde{\Omega}_{1}}\frac{|f(\lambda^{k}x)|^{2}}{\lambda^{2k(\alpha-2)}}dx
=\frac{\frac 1{|B_{\lambda^k}|}\displaystyle\int_{\Omega_{\lambda^k}}|f|^{2}dx}{\lambda^{2k(\alpha-2)}}
\leq[f]^{2}_{\mathcal{L}^{2,-2+\alpha}(0)}\leq(\delta^{*})^{2},
$$
and
$$\|\widetilde{g}\|_{L^{\infty}(\partial_w\widetilde{\Omega}_1)}
=\left\|\frac{g(\lambda^{k}x)}{\lambda^{k\alpha}}\right\|_{L^{\infty}(\partial_w\widetilde{\Omega}_1)}
\leq [g]_{C^\alpha(0)}\leq \delta^{*}.$$
Moreover, since $A$ is $(\delta,1)$-vanishing, we have 
\[
\frac1{|B_1|}\int_{\widetilde\Omega_1}
\left|\widetilde A-(\widetilde A)_{\widetilde\Omega_1}\right|^2\,dx\leq 
\frac1{|B_{\lambda^k}|}\int_{\Omega_{\lambda^k}}
\left|A-(A)_{\Omega_{\lambda^k}}\right|^2\,dx
\leq\delta^2<(\delta^*)^2.
\]
Hence we can apply Lemma $\ref{key22}$ for $w$ to obtain that
$$\left(\frac{1}{|B_{\lambda}|}\int_{\widetilde{\Omega}_{\lambda}}w^{2}dx\right)^{\frac{1}{2}}\leq\lambda^{\alpha}.
$$
We scale back to get
$$\left(\frac{1}{|B_{\lambda^{k+1}}|}\int_{{\Omega}_{\lambda^{k+1}}}|u(x)|^{2}dx\right)^{\frac{1}{2}}
\leq\lambda^{(k+1)\alpha} .
$$
Thus we prove the $(k+1)$-th step.

Finally, for $0<r\leq 1$, there is a k with $\lambda^{k+1}\leq r\leq \lambda^{k}$, and that
\begin{eqnarray*}
\frac{1}{|B_{r}|}\int_{{\Omega}_{r}}|u(x)|^{2}dx&\leq&\frac{1}{\lambda^{n}|B_{\lambda^{k}}|}\int_{{\Omega}_{\lambda^{k}}}|u(x)|^{2}dx\\
&\leq& \frac{\lambda^{2k\alpha}}{\lambda^{n}}\\
&\leq& \frac {1}{\lambda^{n+2\alpha }}r^{2\alpha}.
\end{eqnarray*}
This implies that $u$ is $C^\alpha$ at $0$ in the $L^2$ sense, and
$$
[u]_{\mathcal L^{2,\alpha}(0)}
\leq C\left(\frac{1}{|B_1|}\int_{\Omega_1}u^2\,dx\right)^{1/2}
+C\left([f]_{\mathcal L^{2,-2+\alpha}(0)}
+\|g\|_{C^\alpha(0)}\right)
$$
for some universal constant $C$.
Thus we complete the proof of Theorem \ref{thm1}.

\subsection{\texorpdfstring{$H^1$}{H1} regularity}
To prove Theorem \ref{thm2}, we extend $g$ to a function in
$H^1(\Omega)\cap C(\overline\Omega)$, solve the equation with this
boundary value, and apply uniqueness. We first prove the boundary-layer
estimate for Reifenberg domain.
\begin{lm}\label{rei} 
There exists $\delta_0=\delta_0(n)>0$ such that the following holds.
Let $\Omega \subset \mathbb{R}^n$ be a bounded $(\delta,R)$-Reifenberg
flat domain with $0<\delta<\delta_0$. For any
$x_0\in\partial\Omega$ and $r>0$, let
$E_d(x_0,r) := \Omega_r(x_0) \cap \{x : d(x,\partial\Omega) \le d\}$. Then
$$ |E_d(x_0,r)| \le C r^{n-1} d^{1-C_1\delta^{2}} $$
where $C$ depends on $n,\delta,R$, and $\Omega$; $C_1$ depends only
on $n$. The flatness parameter is small enough that $C_1\delta^{2}<1$.
\end{lm}
\begin{proof}
Without loss of generality, we may assume $x_0=0$. We consider three cases.

\vspace{0.5em}
\noindent \textbf{Case 1: $r \le R/4$ and $d < \beta r$}

Here $\beta\in(0,0.001)$ is a small universal constant. Set
\[
\widetilde\Omega=(2r)^{-1}\Omega,\qquad
\widetilde d=\frac{d}{2r}<\frac{\beta}{2}.
\]
The domain $\widetilde\Omega$ is Reifenberg flat up to scale
$R/(2r)\geq2$.

By Lemma 2.1 of \cite{H14} with $m=n-1$, $P=0$, $R=2$, and
$\varepsilon=\widetilde d$, there exists a smooth $(n-1)$-dimensional
surface
$\widetilde S_{\widetilde d}\subset B_2(0)\cap
N(\partial\widetilde\Omega,\widetilde d)$ such that:
\begin{enumerate}
    \item $\partial\widetilde{\Omega} \cap B_1(0) \subset N(\widetilde S_{\widetilde d},\widetilde d)$;
    \item $\mathcal{H}^{n-1}(\widetilde S_{\widetilde d}) \le C_2 \widetilde d^{-C_1\delta^{2}}$;
    \item $\mathcal{H}^{n-1}(\widetilde S_{\widetilde d} \cap B_{\widetilde d}(y)) \ge C_3 \widetilde d^{n-1}$ for every $y \in \widetilde S_{\widetilde d} \cap B_{1.2}(0)$.
\end{enumerate}

We choose a maximal disjoint subfamily of balls
$\{B_{\widetilde d}(y_i)\}_{i=1}^N$ with centers
$y_i \in \widetilde S_{\widetilde d} \cap B_{1.2}(0)$. By maximality,
these balls enlarged by a factor of $2$ cover the surface:
$$ \widetilde S_{\widetilde d} \cap B_{1.2}(0)
\subset \bigcup_{i=1}^N B_{2\widetilde d}(y_i). $$

By (2) and (3),
$$ N \cdot C_3 \widetilde d^{n-1}
\le \sum_{i=1}^N \mathcal{H}^{n-1}(\widetilde S_{\widetilde d}
\cap B_{\widetilde d}(y_i))
\le \mathcal{H}^{n-1}(\widetilde S_{\widetilde d})
\le C_2 \widetilde d^{-C_1\delta^{2}}. $$
Hence $N \le C \widetilde d^{-(n-1+C_1\delta^{2})}$.

For any point $z \in \partial\widetilde{\Omega} \cap B_1(0)$, property
(1) gives $y \in \widetilde S_{\widetilde d}$ such that
$|z-y|<\widetilde d$. Since $|y|<1+\widetilde d\leq1.2$, the point $y$
belongs to some $B_{2\widetilde d}(y_i)$. Hence
$|z-y_i|<3\widetilde d$, and therefore
$$ \partial\widetilde{\Omega} \cap B_1(0)
\subset \bigcup_{i=1}^N B_{3\widetilde d}(y_i). $$

If $x\in\widetilde\Omega\cap B_{1/2}$ and
$\operatorname{dist}(x,\partial\widetilde\Omega)\leq\widetilde d$,
a nearest boundary point $z$ satisfies $|z|<1/2+\widetilde d<1$.
Thus the preceding boundary covering implies
$x\in B_{4\widetilde d}(y_i)$ for some $i$. Writing
\[
\widetilde E_{\widetilde d}
:=\widetilde\Omega\cap B_{1/2}\cap
\{x:\operatorname{dist}(x,\partial\widetilde\Omega)\leq\widetilde d\},
\]
we obtain
\[
|\widetilde E_{\widetilde d}|
\leq\sum_{i=1}^N\omega_n(4\widetilde d)^n
\leq C\widetilde d^{1-C_1\delta^{2}}.
\]
Scaling back, we obtain
\begin{equation}\label{eq:rei-strong-local}
|E_d(0,r)|=(2r)^n|\widetilde E_{\widetilde d}|
\leq C r^{n-1+C_1\delta^{2}}d^{1-C_1\delta^{2}}.
\end{equation}
Since $r\leq R/4$, the factor $r^{C_1\delta^{2}}$ can be absorbed
into a constant depending on $R$, proving the asserted bound in this case.

\vspace{0.5em}
\noindent \textbf{Case 2: $r \le R/4$ and $d \ge \beta r$} 

Since $E_d(0,r)\subset B_r$,
\[
|E_d(0,r)|\leq\omega_n r^n
\leq\omega_n\beta^{-(1-C_1\delta^{2})}
 r^{n-1+C_1\delta^{2}}d^{1-C_1\delta^{2}}.
\]
This again proves both \eqref{eq:rei-strong-local} and the asserted
bound, since $r\leq R/4$.

\vspace{0.5em}
\noindent \textbf{Case 3: $r > R/4$}\par\nopagebreak[4]

Set $a=R/4$. By compactness, choose a finite covering
\[
\partial\Omega\subset\bigcup_{j=1}^M B_{a/2}(z_j),
\qquad z_j\in\partial\Omega,
\]
where $M$ may depend on $\Omega$ and $R$, but is independent of
$r$ and $d$. If $0<d\leq a/2$, every point of $E_d(0,r)$ lies in
some $E_d(z_j,a)$: choose a nearest boundary point and then a ball
of the covering containing that point. Cases 1 and 2 therefore give
\[
|E_d(0,r)|\leq\sum_{j=1}^M|E_d(z_j,a)|
\leq CM a^{n-1}d^{1-C_1\delta^{2}}
\leq CM r^{n-1}d^{1-C_1\delta^{2}},
\]
because $r>a$. If $d>a/2$, then
\[
|E_d(0,r)|\leq|\Omega|
\leq |\Omega|a^{-(n-1)}(2/a)^{1-C_1\delta^{2}}
 r^{n-1}d^{1-C_1\delta^{2}}.
\]
Together with Cases 1 and 2, this proves the lemma.
\end{proof}

\begin{proof}[Proof of Theorem \ref{thm2}]
Let $\{Q_k\}$ be a Whitney decomposition of $\Omega$, and write
\[
d_k=\operatorname{diam}Q_k,\qquad
\widetilde Q_k=\frac65 Q_k,\qquad
d(x)=\operatorname{dist}(x,\partial\Omega).
\]
Throughout the proof, $C$ is independent of $u,f,g$ and of the
Whitney cubes.

Choose a smooth partition of unity $\{\varphi_k\}$ subordinate to
$\{\widetilde Q_k\}$ such that
\[
0\leq\varphi_k\leq1,\qquad
\sum_k\varphi_k=1\quad\text{in }\Omega,\qquad
\operatorname{supp}\varphi_k\subset\widetilde Q_k,\qquad
|\nabla\varphi_k|\leq Cd_k^{-1}.
\]
The sum is locally finite; see \cite[Chapter VI, Section 1]{Stein70}.
Let $x_k$ be the center of $Q_k$, and choose $y_k\in\partial\Omega$
so that
\[
|x_k-y_k|=\operatorname{dist}(x_k,\partial\Omega).
\]
Define
\[
G(x):=\sum_k\varphi_k(x)g(y_k),\qquad x\in\Omega.
\]
Then $G\in C^\infty(\Omega)$ and
$|G|\leq\|g\|_{L^\infty(\partial\Omega)}$. For
$x\in\widetilde Q_k$, the Whitney properties give
\[
C^{-1}d(x)\leq d_k\leq Cd(x),\qquad |x-y_k|\leq Cd(x).
\]
Hence, for every $y\in\partial\Omega$,
\[
\begin{aligned}
|G(x)-g(y)|
&\leq\sum_k\varphi_k(x)|g(y_k)-g(y)|\\
&\leq C\|g\|_{C^\alpha(\partial\Omega)}|x-y|^\alpha,
\end{aligned}
\]
where we used $d(x)\leq|x-y|$. Thus $G$ extends continuously to
$\overline\Omega$ with $G=g$ on $\partial\Omega$.

Fix $x\in\Omega$ and choose $k_0$ with $\varphi_{k_0}(x)>0$.
Local finiteness gives $\sum_k\nabla\varphi_k(x)=0$, and therefore
\[
\nabla G(x)=\sum_k\nabla\varphi_k(x)
\bigl(g(y_k)-g(y_{k_0})\bigr).
\]
If $\nabla\varphi_k(x)\neq0$, then $|y_k-y_{k_0}|\leq Cd(x)$.
Since the cubes $\widetilde Q_k$ have bounded overlap,
\begin{equation}\label{eq:energy-lifting-gradient}
|\nabla G(x)|
\leq C\|g\|_{C^\alpha(\partial\Omega)}d(x)^{\alpha-1}.
\end{equation}

Fix $x_0\in\partial\Omega$ and take
$r=2\operatorname{diam}(\Omega)$ in Lemma \ref{rei}. Since
$\Omega\subset B_r(x_0)$, that lemma gives
\[
|\{x\in\Omega:d(x)\leq t\}|
\leq Ct^{1-C_1\delta^{2}},\qquad t>0,
\]
where $C$ may depend on $\Omega$. Set
$t_j=2^{-j}$ and $E_j=\{x\in\Omega:t_{j+1}<d(x)\leq t_j\}$.
Since $\alpha<1$, we obtain
\[
\begin{aligned}
\int_{\{d(x)\leq1\}}d(x)^{2\alpha-2}\,dx
&\leq C\sum_{j=0}^{\infty}t_j^{2\alpha-2}|E_j|\\
&\leq C\sum_{j=0}^{\infty}
2^{-j(2\alpha-1-C_1\delta^{2})}<\infty.
\end{aligned}
\]
The series converges because $\delta<\delta^{**}$.
On $\{d(x)>1\}$, the integral is bounded by $|\Omega|$.
Together with \eqref{eq:energy-lifting-gradient} and the bound for
$G$, this proves
\begin{equation}\label{eq:energy-lifting-bound}
G\in H^1(\Omega)\cap C(\overline\Omega),\qquad
\|G\|_{H^1(\Omega)}
\leq C\|g\|_{C^\alpha(\partial\Omega)}.
\end{equation}

By Lax--Milgram Theorem, there exists $w\in H_0^1(\Omega)$ such that
\[
\int_\Omega a_{ij}D_jwD_i\psi\,dx
=\int_\Omega f\psi\,dx
-\int_\Omega a_{ij}D_jGD_i\psi\,dx
\qquad\forall\,\psi\in H_0^1(\Omega).
\]
Testing with $w$, and using ellipticity and Poincar\'e's inequality,
we obtain
\[
\|w\|_{H^1(\Omega)}
\leq C\bigl(\|f\|_{L^2(\Omega)}+\|\nabla G\|_{L^2(\Omega)}\bigr).
\]
Set $U=G+w$. Then we have
\begin{equation}\label{eq:energy-variational-equation}
\int_\Omega a_{ij}D_jUD_i\psi\,dx
=\int_\Omega f\psi\,dx
\qquad\forall\,\psi\in H_0^1(\Omega),
\end{equation}
and 
\begin{equation}\label{eq:energy-variational-bound}
\|U\|_{H^1(\Omega)}
\leq C\left(\|f\|_{L^2(\Omega)}
+\|g\|_{C^\alpha(\partial\Omega)}\right).
\end{equation}

The boundary truncation argument in the proof of Theorem
\ref{thm:existence} applies to $U-G\in H_0^1(\Omega)$ with
$G\in H^1(\Omega)\cap C(\overline\Omega)$ and shows that $U$
satisfies Definition \ref{uni}. Since $\Omega$ is $\theta$-admissible,
Theorem \ref{thm:uniqueness} gives $u=U$ a.e. in $\Omega$.
The conclusion follows from \eqref{eq:energy-variational-bound}.

\end{proof}

\begin{rem}[Global $W^{1,p}$ regularity]\label{rem:global-w1p}
Let $2\leq p<\infty$ and $0<\alpha<1$. There exists
$\delta_p=\delta_p(n,\Lambda,p)>0$ such that the following holds.
Let $\Omega$ be a bounded $(\delta,1)$-Reifenberg flat domain, where
\[
0<\delta<\min\{\delta_0,\delta_p\},\qquad
p(1-\alpha)<1-C_1\delta^{2}.
\]
Here $C_1=C_1(n)$ and $\delta_0=\delta_0(n)>0$ are as in Lemma
\ref{rei}. Suppose that $A$ satisfies \eqref{ellipticity} and is
$(\delta,1)$-vanishing, $f\in L^p(\Omega)$, and
$g\in C^\alpha(\partial\Omega)$. If $u$ is a weak solution of
\eqref{1} in the sense of Definition \ref{uni}, then
$u\in W^{1,p}(\Omega)$ and
\begin{equation}\label{eq:global-w1p-estimate}
\|u\|_{W^{1,p}(\Omega)}
\leq C\left(\|f\|_{L^p(\Omega)}
+\|g\|_{C^\alpha(\partial\Omega)}\right),
\end{equation}
where $C$ depends only on $n,\Lambda,p,\alpha,\delta$, and $\Omega$.

Indeed, the function $G$ constructed in the proof of Theorem
\ref{thm2} belongs to $C(\overline\Omega)$ and equals $g$ on
$\partial\Omega$. By \eqref{eq:energy-lifting-gradient} and Lemma
\ref{rei},
\[
\int_{\{d(x)\leq1\}}|\nabla G|^p\,dx
\leq C\|g\|_{C^\alpha(\partial\Omega)}^p
\sum_{j=0}^{\infty}
2^{-j(1-C_1\delta^{2}-p(1-\alpha))}<\infty.
\]
The remaining part of $\Omega$ is controlled by
\eqref{eq:energy-lifting-gradient}. Thus
\begin{equation}\label{eq:wp-whitney-bound}
\|G\|_{W^{1,p}(\Omega)}\leq C\|g\|_{C^\alpha(\partial\Omega)}.
\end{equation}

Since $p\geq2$ and $\Omega$ is bounded, we have
$f\in L^2(\Omega)$ and $G\in H^1(\Omega)\cap C(\overline\Omega)$.
The proof of Theorem \ref{thm2} therefore gives $U=G+w$, with
$w\in H_0^1(\Omega)$, satisfying Definition \ref{uni}.
By Theorem \ref{thm:uniqueness}, $u=U$.

Choose a ball $D\supset\overline\Omega$, extend $f$ by zero to $D$,
and let $z\in W^{2,p}(D)\cap W_0^{1,p}(D)$ solve
$-\Delta z=\widetilde f$ in $D$. The standard estimate on a ball gives
\[
F:=-\nabla z,\qquad
\operatorname{div}F=f\quad\text{in }\Omega,\qquad
\|F\|_{L^p(\Omega)}\leq C\|f\|_{L^p(\Omega)}.
\]
Consequently,
\[
-D_i(a_{ij}D_jw)=D_i(F_i+a_{ij}D_jG)
\quad\text{in }\Omega,\qquad w\in H_0^1(\Omega).
\]
For $p>2$, \cite[Theorem 1.5]{BW04} and uniqueness of the energy
solution give $w\in W_0^{1,p}(\Omega)$ and
\[
\|\nabla w\|_{L^p(\Omega)}
\leq C\left(\|F\|_{L^p(\Omega)}
+\|\nabla G\|_{L^p(\Omega)}\right).
\]
That theorem allows nonsymmetric matrices
\cite[Remark 2.2]{BW04}. We choose $\delta_p$ sufficiently small
for that theorem and for the exterior measure-density condition.
For $p=2$, the same estimate follows directly by testing the
variational equation with $w$. Poincar\'e's inequality,
$u=G+w$, and \eqref{eq:wp-whitney-bound} prove
\eqref{eq:global-w1p-estimate}. In particular, for $p=2$ the
small-BMO assumption is unnecessary, as already shown in Theorem
\ref{thm2}.

The lower threshold $\alpha>1-1/p$ is sharp even for the
Laplacian on smooth bounded domains in $\mathbb R^n$;
see Example \ref{ex:sharp-w1p}.

\end{rem}

\begin{ex}[Sharpness of the lower H\"older thresholds]\label{ex:sharp-w1p}
The counterexample in \cite[Proposition 11]{BDL26} extends to $n\geq3$
by making it independent of the remaining coordinates.
Choose a smooth bounded domain in $\mathbb R^2\times\mathbb R^{n-2}$ with
\[
B_1^2\times B_1^{n-2}\subset\Omega\subset B_1^2\times B_2^{n-2},
\]
where the superscripts denote dimensions. For $1<p<\infty$ and
$\alpha=1-1/p$, set
\[
u(r\cos\theta,r\sin\theta,y)
=\sum_{j=1}^{\infty}2^{-j\alpha}r^{2^j}\cos(2^j\theta).
\]
By that proposition and \cite[Section 5]{Zhu07}, respectively,
\[
\Delta u=0,\qquad u\in C^\alpha(\overline\Omega),\qquad
\int_\Omega|\nabla u|^p\,dx
\geq\int_{B_1^2\times B_1^{n-2}}|\nabla u|^p\,dx=\infty,
\]
since $\sum_j2^{-j\alpha p}2^{j(p-1)}=\sum_j1=\infty$.
For $g=u|_{\partial\Omega}$, uniform approximation by the partial sums
verifies Definition \ref{uni} as in the proof of Theorem
\ref{thm:existence}. Dilating $\Omega$ gives
arbitrarily small Reifenberg flatness at scales at most one.
Taking $p=2$ proves sharpness of the threshold in Theorem \ref{thm2}.

\end{ex}

\end{document}